\documentclass{article}
\usepackage{graphicx}        %
\usepackage{multicol}        %
\usepackage{multirow}        %
\usepackage{array}
\usepackage{float}
\newcolumntype{P}[1]{>{\centering\arraybackslash}p{#1}}

\usepackage{pgfplotstable}
\usepackage{pgfplots}
\pgfplotsset{compat=1.16}

\usepackage{booktabs}
\usepackage{tikz}
\usetikzlibrary{patterns}

\usepackage{amsmath,amsfonts}
\usepackage{nicefrac}

\newcommand{\R}{\mathbb{R}}
\newcommand{\Ogrande}{\mathcal{O}}

\newcommand{\OSC}{\mathrm{I}}
\newcommand{\RBF}{\ensuremath{\mathsf{RBF}}}
\newcommand{\WENO}{\ensuremath{\mathsf{WENO}}}
\newcommand{\Qone}{\ensuremath{\mathsf{Q1}}}

\newcommand{\pder}[2]{{#1}_{#2}}

\usepackage{algorithm}

\usepackage{url}

\begin{document}

\title{Surface reconstruction from point cloud using a semi-Lagrangian scheme with local interpolator}
\author{Silvia Preda%
 \thanks{Dipartimento di Scienza e Alta Tecnologia, Universit\`a dell'Insubria, Via Valleggio 11, 22100 Como (Italy)
   {\tt e-mail: spreda@uninsubria.it} - ORCID 0009-0003-8405-2245}
 \and Matteo Semplice%
 \thanks{Dipartimento di Scienza e Alta Tecnologia, Universit\`a dell'Insubria, Via Valleggio 11, 22100 Como (Italy)
   {\tt e-mail: matteo.semplice@uninsubria.it} - ORCID 0000-0002-2398-0828}}
\maketitle

\begin{abstract}
We propose a level set method to reconstruct unknown surfaces from point clouds, without assuming that the connections between points are known. We consider a variational formulation with a curvature constraint that minimizes the surface area weighted by the distance of the surface from the point cloud. More precisely we solve an equivalent advection-diffusion equation that governs the evolution of an initial surface described implicitly by a level set function. Among all the possible representations, we aim to compute the signed distance function at least in the vicinity of the reconstructed surface. The numerical method for the approximation of the solution is based on a semi-Lagrangian scheme whose main novelty consists in its coupling with a local interpolator instead of a global one, with the aim of saving computational costs. In particular, we resort to a multi-linear interpolator and to a Weighted Essentially Non-oscillatory one, to improve the accuracy of the reconstruction.
Special attention has been paid to the localization of the method and to the development of fast algorithms that run in parallel, resulting in faster reconstruction and thus the opportunity to easily improve the resolution. A preprocessing of the point cloud data is also proposed to set the parameters of the method. Numerical tests in two and three dimensions are presented to evaluate the quality of the approximated solution and the efficiency of the algorithm in terms of computational time.

{{\bf Keywords.} Curvature regularized surface reconstruction
\and Point cloud
\and Level set methods
\and Semi-Lagrangian schemes
\and $\WENO$ interpolation
}
\end{abstract}

\section{Introduction}
\label{sec:intro}
The problem of acquiring, creating and processing 3d digital models of real objects has gained importance in numerous applications across diverse industries and areas of research including physics, computer graphics, medical imaging, urban planning and cultural heritage. Especially in the latter, real artistic manufacts are not Cartesian domains, instead they usually present complicated geometries and topologies and they are not built from precise (e.g. CAD) drawings. Therefore, their shape must be acquired from the real object itself, e.g. via 3d laser scanning or photogrammetry \cite{scan3d2020,Re:2011,ReEl:2006}.

In many cases, the only information available for the reconstruction of the original shape is a set of unorganized and non-uniformly spaced points that does not contain any information about the ordering or the connection between them. This makes surface reconstruction from point clouds complex and time-consuming. Additionally, as an infinite number of surfaces may pass through or near the data points, this problem turns out to be ill-posed, with no unique solution. The challenge is thus to get a good approximation of the data set that should have some smoothness properties while being able to retain as many details and features of the real object as possible. Moreover, the reconstruction should be also useful for dynamic operations and not only for static representations.

A prototypical application of this topic in the field of cultural heritage has been presented in \cite{cdss:mach19}, where a computational domain having the shape of a work of art has been computed from a set of unorganized points and then employed for the simulation of a reaction-diffusion model investigating the evolution of damage caused by the interaction between pollutants in the surrounding environment with the work of art itself. In fact, mathematical models for marble sulfation or other reaction-diffusion mechanisms of degradation of cultural heritage \cite{CSS:monum,cdss:mach19} require to perform computations on complex domains representing a work of art and this is where level set functions find an interesting application.
In particular, our main interest is to produce high quality level set descriptions of complex objects, which can be later used as domain definition in partial differential equations (PDEs) computations via ghost-cell methods \cite{OshStFe:2004,GF:2002,CR2018}.

In general, there are two approaches for representing a surface: the explicit one and the implicit one. The former include mesh-based techniques, such as Delaunay triagulations and Voronoi diagrams \cite{Ed:delaunay1998,Am:delaunay2004,CaFr:surveyDelaunay} and parametric techniques, e.g. NURBS \cite{NURBS1,NURBS2}. Although they represent some of the most explored techniques in this topic, the main problem with this approach is to ensure a watertight reconstruction in cases when an evolution of the surface occurs. Recent works, for example \cite{DaMi:surfMesh2015}, propose to overcome this difficulty by considering an additional mesh adjustment in order to control the quality of the mesh during the evolution. However, the state of the art agrees that explicit reconstructions are difficult to handle when dealing with moving boundaries or complex and changing topologies. 

On the other hand, using an implicit approach, one has the advantage of better handling topological flexibility, while having a very simple data structure that also allows very simple Boolean operations on the detected surface. Among implicit methods, traditional ones enclose interpolation techniques to approximate the desired implicit function as a combination of some smooth basis functions. The problem here is obviously related to the large linear systems that need to be solved, resulting in high computational costs and, even worse, in ill-conditioned matrices when the number of interpolated points is very large. Radial Basis Functions ($\RBF$) with global and local support have found many applications along this line \cite{RBF1,RBF2,RBF3} and some least-square-based methods \cite{leastSquare2} belong to this family too. Other methods use local estimators to associate an oriented plane to each point in the cloud and thus approximate the signed distance function representing the surface \cite{hoppe1992}. A different and widespread approach for the processing of point cloud data is constituted by the application of level set methods \cite{OshSet:1988}. This last line of research constitutes the focus of this work. 

For a broader overview about the reconstruction of real objects starting from point clouds, the reader can refer to \cite{BeTaSe:survey2014,BeTaSe:survey2016}. In these surveys various algorithms have been compared and classified according to their ability to handle imperfections in the point clouds, their input requirements, the class of shape that they are able to recover and the type and the quality of the final reconstruction. Significant advancements in surface reconstruction are also being driven by the integration of machine learning techniques, especially deep learning models, which have shown remarkable results in learning implicit surface representations directly from unorganized point clouds \cite{DL:2021,DL:2023gen,DL:2023apr}. 

Since their introduction, level set methods have emerged as a powerful and versatile tool in a wide range of applications \cite{SSO:1994,Se:1999,OshStFe:2004} including image processing and surface reconstruction \cite{Zhao:2000,LSApp:2004}. The central idea is to represent both an $n$-dimensional object $\Xi$ and its boundary by a so called level set function $\phi:\R^n\rightarrow\R$ such that $\Xi = \{ \vec{x}\in\R^n : \phi(\vec{x})<0\}$ and its boundary $\partial\Xi$ is the zero isocontour of $\phi$. $\Xi$ and $\partial\Xi$ can be then evolved in time using PDEs for $\phi$, driving the deformation of an initial surface. With this approach, one can easily handle noisy and incomplete point clouds, since the level set function evolves smoothly and naturally adapts to the input data, enjoying all the geometrical and topological advantages of working with implicit representations.

Our work delves into the approximation of the solution of the level set evolution PDE presented in \cite{Zhao:2000} that leads the evolution of an initial guess to get the final, steady-state shape. The model is based on the minimization of an energy functional penalized by a curvature term. This allows us to trade off the exact vanishing of the level set function on the points of the cloud for a control on the maximum curvature of the resulting zero level set surface. This curvature term is also expected to be especially useful in counteracting noise in the data.  For the numerical approximation, we apply a semi-Lagrangian scheme \cite{CaFe:2017}, first presented in \cite{FaFe:2003} for curvature-related equations, with the application of a local interpolator for space reconstructions. Here we employ both a multilinear interpolation and a third order accurate Weighted Essentially Non-oscillatory ($\WENO$) reconstruction \cite{CFR05}. We find that both are faster than the $\RBF$ approach of \cite{CaFe:2017}, with the second one yelding more satisfactory results.
Since it is preferable to have a level set function that is close to the signed distance function, we also perform a constrained reinitialization \cite{HaMeSc:2008} to keep the gradient of the solution well behaved and to allow reliable computation of basic geometrical quantities such as normals, curvature and surface area.
Finally, since we deal with three spatial dimensions, and keeping in mind the application to complex shapes occurring in practice, we will employ distributed memory parallelization and some other strategies to reduce the computational cost of our scheme.

The paper is organized as follows. In \S\ref{sec:mathmodel} we briefly recall the mathematical model \cite{Zhao:2000} and the level set formulation \cite{OshSet:1988} that lay behind our governing PDE. In \S\ref{sec:semilagrangian} all the numerical schemes used in the algorithm are detailed, together with our strategies to save computational costs. In \S\ref{sec:numtests} two- and three-dimensional tests of increasing complexity are presented. Finally \S\ref{sec:conclusions} draws some conclusions and highlights future perspectives.

\section{Mathematical model}
\label{sec:mathmodel}
\newcommand{\pcloud}{\mathcal{S}}

Let us assume that we are given a set of points $\pcloud=\{{Q}_1,\ldots,{Q}_N\}$ in a bounded region of $\R^n$ 
and let us define $d(\vec{x}) = \min_{Q \in \pcloud}|\vec{x}-Q |$ to be the distance function to $\pcloud$. Our model for surface reconstruction from point clouds is based on the minimization of the surface energy
\begin{equation}\label{eq:energy}
  E_p(\Gamma) = \Big( \int_{\Gamma} d^p(\vec{x}) ds \Big) ^{1/p},   
\end{equation}
where $\Gamma$ is a closed surface of co-dimension one in $\R^n$. To achieve the minimum and get the final shape, we continuously deform an initial surface $\Gamma_0$ following the gradient descent of the energy functional. In this process, the evolving surface $\Gamma(t)$ is represented implicitly using a level set function to capture the moving interface (see \cite{Zhao:2000} for full details).
\newline

Let $\Xi(t)$ be the region enclosed by $\Gamma(t)$. Level set methods, first introduced by Osher and Sethian \cite{OshSet:1988}, consider a level set function $\phi(\vec{x},t)$ associated with $\Xi(t)$ such that
\begin{equation}
    \begin{split}
        \phi(\vec{x},t) &<0 \quad \text{inside} \quad \Xi(t),\\
        \phi(\vec{x},t) &=0 \quad \text{on} \quad \Gamma(t),\\
        \phi(\vec{x},t) &>0 \quad \text{outside} \quad \Xi(t).
    \end{split} 
\end{equation}    

In terms of the level set function, the energy functional \eqref{eq:energy} can be rewritten as
\begin{equation} \label{eq:energy:ls}
    E_p(\phi) = \Bigg( \int_{\Xi} \lvert d(\vec{x}) \rvert ^p \delta (\phi) \lvert \nabla \phi \rvert d\vec{x} \Bigg)^{1/p}
\end{equation}
where $\delta$ denotes the Dirac-delta function.
Following \cite{Zhao:2000}, we minimize the energy \eqref{eq:energy:ls} by computing the solution of the following
evolutionary problem
\begin{equation}\label{eq:ZhaoPDE}\small
    \begin{split}
        \pder{\phi}{t}(\vec{x},t) 
        &= \Bigg[  \frac{d(\vec{x})}{E_p(\phi)}   \Bigg]^{p-1} \Bigg( \nabla d(\vec{x})\cdot \nabla \phi(\vec{x},t) + \frac{1}{p} d(\vec{x}) \nabla\cdot\left( \frac{\nabla \phi(\vec{x},t)}{|\nabla \phi(\vec{x},t)|} \right)|\nabla \phi(\vec{x},t)| 
          \Bigg),\\
        \phi(\vec{x},0) &= \phi_0(\vec{x}), 
    \end{split}
\end{equation}
where $\phi_0(\vec{x})$ is a suitable initial data such that $\Gamma_0=\{\vec{x}\in\R^n: \phi_0(\vec{x})=0\}$.

In the above equation, the term $\nabla d(\vec{x})\cdot \nabla \phi(\vec{x},t)$ drives the surface towards the dataset $\pcloud$, while the second term tempers the maximal curvature of $\Gamma$. Changing the balance between them can lead to surfaces closer to the point cloud but with sharper edges or to more rounded surfaces that are a little further away from $\pcloud$.

The parameter $p$ controls the relative influence of the curvature term and of the global scaling factor $\big[ \frac{d(\vec{x})}{E_p(\phi)} \big]^{p-1}$ which makes the surface to move quicker where it is further away from the point cloud and slower when it approaches $\pcloud$. We point out that the $E_p(\phi)$ denominator slows down the entire evolution when the functional is large; this can be an issue when the initial surface is very far from the cloud.

In \cite{Kosa2017}, a model similar to \eqref{eq:ZhaoPDE} is employed, wherein the factor $d(\vec{x})$ in front of the curvature term is replaced by a constant $\delta\in[0,1]$ to control the balance between the two terms. 
While we prefer to keep the $d(\vec{x})$ factor in the evolution equation which will slow down the evolution in the vicinity of the point cloud, we also introduce an additional parameter $\delta$ multiplying $d(\vec{x})$ in the curvature term. Our final formulation then reads
\begin{equation}\label{eq:levelset:pde}\small
    \begin{split}
        \pder{\phi}{t}(\vec{x},t) 
        &= \Bigg[  \frac{d(\vec{x})}{E_p(\phi)}   \Bigg]^{p-1} \Bigg( \nabla d(\vec{x})\cdot \nabla \phi(\vec{x},t) + \frac{\delta}{p} d(\vec{x}) \nabla\cdot\left( \frac{\nabla \phi(\vec{x},t)}{|\nabla \phi(\vec{x},t)|} \right)|\nabla \phi(\vec{x},t)| \Bigg),\\
        \phi(\vec{x},0) &= \phi_0(\vec{x}).
    \end{split}
\end{equation}
In this way our numerical algorithm can be used in both regimes.
If we set $\delta=0$, the model becomes purely convective and leads quickly to a reconstructed surface close to a piecewise linear approximation. In constrast, the case with $\delta=1$, contains a weighted curvature regularization effect, is more computationally expensive, but leads to a smoother reconstructed surface.

Finally,
we observe that for visualization purposes it is often sufficient to compute a function $\phi$ such that its zero level set is close to the data set $\pcloud$. Equation \eqref{eq:levelset:pde} is apt for this purpose, but its transport term tends to produce high gradients in the computed function $\phi$. However, for computing normals of the reconstructed surface or to use the level set function as a mathematical description of a domain into which a PDE solver has to be applied, it is more robust to compute a level set function with moderate gradients. In particular, a level set function with $|\nabla\phi|=1$ is called a signed-distance function and our algorithm will be designed to ensure that the computed solution is a signed distance at least in the vicinity of the point cloud.

\section{Numerical scheme}
\label{sec:semilagrangian}
\newcommand{\grid}{\mathcal{G}}
\newcommand{\dt}{\mathrm{\Delta}t}
\newcommand{\dx}{\mathrm{\Delta}x}
\newcommand{\ddelta}{\mathrm{\Delta}}
\newcommand{\cupdot}{\mathbin{\mathaccent\cdot\cup}}
The numerical evolution of \eqref{eq:levelset:pde} is computed following the semi-Lagrangian approach presented in \cite{CaFe:2017}. As such, the scheme is explicit, but not constrained by a parabolic-type CFL condition.
However, the interpolation operator proposed in \cite{CaFe:2017} involves the solution of a large linear system at each timestep. For this reason, we will propose the use of a local interpolant, which is more efficient, especially in view of a parallel implementation.

The scheme considers a background Cartesian mesh $\grid$ with uniform mesh width $\dx$  on a bounded region of $\R^2$ or $\R^3$, containing the point cloud $\pcloud$. Let $\vec{x}_j\in\grid$ be a point in the Cartesian mesh and let us denote by $\phi^n_j$ the approximate value of $\phi$ at $\vec{x}_j$ and time $t^n$. 

\begin{algorithm}
\caption{Computing the levelset $\phi$ from a point cloud $\pcloud$}\label{algo}
\vspace{5pt}
\begin{enumerate}
\item Create a Cartesian grid $\grid$ encompassing all of $\pcloud$
\item \label{algo:d} Compute the distance function $d(\vec{x})$ at all grid points, see \S\ref{ssec:distance}
\item \label{algo:u0} Set initial data $\{\phi^0_j\}_{\vec{x}_j\in\grid}$ as in \S\ref{subsec:initialData}

\item Loop:
  \begin{enumerate}
    \item \label{algo:banda} Choose computational subgrid $\widetilde\grid\subset\grid$ as in \S\ref{ssec:mask}
    \item \label{algo:energy} Compute the energy functional \eqref{eq:energy:ls} with one of the methods of \S\ref{ssec:energy}
    \item \label{algo:sl} Compute  $\{\phi^{n+1}_j\}_{\vec{x}_j\in\widetilde\grid}$ using \eqref{eq:SL2d} or \eqref{eq:SL3d} from \S\ref{ssec:semilagrangian} and the interpolator \eqref{eq:IQ1} or \eqref{IWENO} from \S\ref{ssec:interp}
    \item Reinitialize $\{\phi^{n+1}_j\}_{j\in\overline\grid}$ as in \S\ref{ssec:reinit} and cut\label{algo:reinit} with \eqref{eq:cut}
  \end{enumerate}
\end{enumerate}
\end{algorithm}

The semi-Lagrangian scheme will be recalled in \S\ref{ssec:semilagrangian} and the interpolation operator will be described in \S\ref{ssec:interp}. The later subsections will describe the other algorithms used in our code: the computation of the distance function $d(\vec{x})$ from the point cloud, the evaluation of the energy functional $E_p(\phi)$, the choice of the initial data $\phi_0$, the reinitialization procedure and the use of cut-off functions to reduce the computational effort. The complete algorithm is summarized in Algorithm~\ref{algo}.
Unless otherwise specified, all numerical gradients needed in the algorithms are computed by centered finite differences.

\subsection{Semi-Lagrangian scheme}\label{ssec:semilagrangian}
The scheme follows a semi-Lagrangian approach both for the advection and for the diffusion terms. In particular, the diffusion term is discretized by averaging the data in a region of size $\sqrt{\dt}$ around the foot of the characteristic of the advection term, as first proposed in \cite{BF:14:SLdiffusion,BCCF:21}. 
\subsubsection*{2d Case}
We compute the update of $\phi^n_j$ as
\begin{equation}\label{eq:SL2d}
    \begin{aligned}
       \phi^{n+1}_j
        \, &=\; \frac12 \sum_{i=1}^2I[\phi^n]\left( \vec{x}^*_{j,i} \right),
        \\
        \vec{x}^*_{j,i} &= \vec{x}_j + C^n_j \dt \nabla d(\vec{x}_j) + \sqrt{\frac{C^n_j \, \delta \, d(\vec{x}_j) \dt}{p}}\,\sigma^n_j \lambda_i,
    \end{aligned}
\end{equation}
where $\lambda_i$ ranges over $\{ -1, +1\}$ and $C^n_j$ is the scale factor $\big[\frac{d(\vec{x_j})}{E_p(\phi^n)}\big]^{p-1}$. The operator $I[\phi^n](\vec{x})$ denotes an interpolation at point $\vec{x}$ of the data $\{\phi^n_j: \vec{x}_j\in\grid\}$. This operator will be specified later. $\sigma^n_j$ denotes the unit vector tangent to the level sets of $\phi$: it is thus orthogonal to the gradient of the level set function and is given by
\begin{equation}\label{sigma2d}
    \sigma^n_j = \frac{1}{\lvert\nabla\phi\rvert} 
    \begin{pmatrix} 
        &\partial_2\phi \\ -&\partial_1\phi
    \end{pmatrix}.
\end{equation}
In equation \eqref{eq:SL2d}, the interpolation points are obtained by adding two terms: the first is the foot of the characteristic pertaining to the advection term in \eqref{eq:levelset:pde}, while the second is a further displacement that generates a diffusion along the tangent space of the level sets, thus discretizing the curvature term in \eqref{eq:levelset:pde}. This discretization of the second order operator has been described in \cite{FaFe:2003,CaFaFe:2010}.

\subsubsection*{3d Case}
Similarly to the previous case, the update of $\phi^n_j$ is computed as
\begin{equation}\label{eq:SL3d}
    \begin{aligned}
       \phi^{n+1}_j
        \, &=\; \frac14 \sum_{i=1}^4I[\phi^n]\left( \vec{x}^*_{j,i} \right),
        \\
        \vec{x}^*_{j,i} &= \vec{x}_j + C^n_j \dt \nabla d(\vec{x}_j) + \sqrt{\frac{C^n_j \, \delta \, d(\vec{x}_j) \dt}{p}}\,\sigma^n_j \lambda_i,
    \end{aligned}
\end{equation}
where again we treat the scale factor $C^n_j$ as a constant.
In the above equation,
$\sigma^n_j = \left[ \nu_1(\nabla \phi^n_j) , \nu_2(\nabla \phi^n_j) \right]
$ is a $3\times2$ matrix whose columns 
span the $2d$ space orthogonal to $\nabla\phi^n$ 
and the row vector $\lambda_i$ ranges over $\{(\pm1,\pm1)\}$,
so that $\sigma^n_j \lambda_i$ represent four points in the local tangent plane.

In order to avoid numerical singularities, one regularizes the
two orthonormal eigenvectors of the projection %
onto the plane tangent to the level sets of $\phi$ as
\begin{align*}
    \nu_1 &= 
    \begin{cases}
    \tilde\nu_1 &\text{if} \, \sqrt{(\partial_1 \phi)^2+(\partial_3 \phi)^2} \neq 0
    \\
    (1,0,0)^T &\text{otherwise}
    \end{cases}
\\    
    \nu_2 &= 
    \begin{cases}
    \tilde\nu_2 &\text{if} \, \sqrt{(\partial_1 \phi)^2+(\partial_3 \phi)^2} \neq 0
    \\
    (0,0,1)^T &\text{otherwise}
    \end{cases}
\end{align*}
with $\tilde\nu_k$ denoting the exact eigenvectors, which are given by
\begin{equation}\label{sigma3d}
\tilde{\nu}_1 = 
\begin{pmatrix}
\tfrac{-\partial_3 \phi}{\sqrt{(\partial_1 \phi)^2+(\partial_3 \phi)^2}}
\\
\\
0
\\
\\
\tfrac{\partial_1 \phi}{\sqrt{(\partial_1 \phi)^2+(\partial_3 \phi)^2}}
\end{pmatrix},
\quad
\tilde{\nu}_2 = 
\frac{1}{|\nabla \phi|}
\begin{pmatrix}
\tfrac{-\partial_1 \phi\,\partial_2 \phi}{\sqrt{(\partial_1 \phi)^2+(\partial_3 \phi)^2}}
\\
\\
\sqrt{(\partial_1 \phi)^2+(\partial_3 \phi)^2}
\\
\\
\tfrac{-\partial_2 \phi\,\partial_3 \phi}{\sqrt{(\partial_1 \phi)^2+(\partial_3 \phi)^2}}
\end{pmatrix}
.
\end{equation}

Because of the factor in front of $\sigma^n_j$ in \eqref{sigma2d}, respectively in front of $\tilde\nu_2$ in \eqref{sigma3d}, the numerical schemes \eqref{eq:SL2d} and \eqref{eq:SL3d} would be singular in cases where $|\nabla \phi^n_j|$ is close to zero. Thus, when $|\nabla \phi^n_j|<C\dt^\alpha$, the schemes are replaced by
\begin{equation}\label{eq:SL:2}
\phi^{n+1}_j = \frac{1}{\lvert \mathcal{N}_j \vert} \sum_{\vec{x}_i\in\mathcal{N}_j} I[\phi^n](\vec{x}_i) 
\end{equation}
where $\mathcal{N}_j$ is the set of the first neighbours of $\vec{x}_j$ in the Cartesian grid $\grid$ and $\lvert \mathcal{N}_j \vert$ represents its cardinality.
In all computations presented in this work, we set $C=10^{-3}$ and $\alpha=1$.

In order to fully specify the numerical algorithm, a suitable interpolation operator has to be specified.

\subsection{Interpolation}\label{ssec:interp}

An important issue for the accuracy and the efficiency of the numerical scheme is clearly represented by the choice of the interpolation operator $I[\phi^n]$. 
In \cite{CaFe:2017}, the authors employ a global $\RBF$ interpolation of the data on $\grid$ defined as
\begin{equation*}
    I_{\text{\RBF}}[\phi^n](\vec{x}) = c_0(\phi^n) + \vec{c}(\phi^n) \cdot \vec{x} + \sum_{\vec{x}_i\in\grid} \mu_i(\phi^n) \psi(|\vec{x}-\vec{x}_i|)
.
\end{equation*}
The coefficients $c_0\in\R,\vec{c}\in\R^n$ and $\mu_i$ of the $\RBF$ interpolator are computed by imposing interpolation conditions at all (or a subset of) points of the grid ($I[\phi](\vec{x}_k)=\phi_k$ for all $\vec{x}_k\in\grid$)
and the additional conditions
$\sum_{\vec{x}_i\in\grid} \mu_i = 0$,
$\sum_{\vec{x}_i\in\grid} x_i \mu_i =0$,
$\sum_{\vec{x}_i\in\grid} y_i \mu_i =0$,
$\sum_{\vec{x}_i\in\grid} z_i \mu_i =0$.
The interpolation is thus computed by solving, at each time step, a linear system with matrix of the form
$\begin{bmatrix} B & P \\ P^T & 0 \end{bmatrix}$, 
where $B$ is an $N\times N$ block ($N=|\grid|$) and $P$ is $N\times 4$.

However, the computation of the $\RBF$ interpolator then becomes a bottleneck of the algorithm due to the solution of the linear system involved. Furthermore, the global linear term $c_0(\phi^n) + \vec{c}(\phi^n) \cdot \vec{x}$, i.e. the blocks $P$ and $P^T$ in the system, forms a strong coupling of all the equations which is difficult to handle for parallel runs as it requires a lot of inter-processor communications.

\begin{figure}
\centering
\begin{tikzpicture}
\usetikzlibrary{patterns}
\definecolor{stencil}{gray}{0.8}
\colorlet{cella}{red}

\begin{scope}
\filldraw[fill=stencil](-0.5,-0.5) rectangle (0.5,0.5);
\fill[pattern=north east lines, pattern color=cella] (-.5,-.5) rectangle (.5,.5);
\foreach \x in {-1.5,-0.5,0.5,1.5}
{
  \draw (\x,-1.75) -- (\x,1.75);
  \draw (-1.75,\x) -- (1.75,\x);
}
\foreach \x in {-0.5,0.5}
{
    \foreach \y in {-0.5,0.5}
    \fill[black] (\x,\y) circle (2pt);
}
\draw [blue, xshift=-.6cm,yshift=-.6cm] 
(-0.15,-0.15) 
-- (-0.15,1.35) 
-- (1.35,1.35) 
-- (1.35,-0.15) 
-- cycle;
\end{scope}
\end{tikzpicture}  \hspace{1cm}
\begin{tikzpicture}
\usetikzlibrary{patterns}
\definecolor{stencil}{gray}{0.8}
\colorlet{cella}{red}

\begin{scope}
\filldraw[fill=stencil](-1.5,-1.5) rectangle (1.5,1.5);
\fill[pattern=north east lines, pattern color=cella] (-.5,-.5) rectangle (.5,.5);
\foreach \x in {-1.5,-0.5,0.5,1.5}
{
  \draw (\x,-1.75) -- (\x,1.75);
  \draw (-1.75,\x) -- (1.75,\x);
}
\foreach \x in {-1.5,-0.5,0.5,1.5}
{
    \foreach \y in {-1.5,-0.5,0.5,1.5}
    \fill[black] (\x,\y) circle (2pt);
}

\draw [orange, xshift=-1.65cm,yshift=-1.65cm] 
(0,0) 
-- (0,3.3) 
-- (3.3,3.3) 
-- (3.3,0) 
-- (0,0);
\end{scope}
\end{tikzpicture} 
\caption{Stencils of the two-dimensional $\Qone$ and $\WENO$ reconstructions. The red hatched region represents the cell $\Omega$ in which we compute the reconstruction. The multi-linear interpolator only requires the vertices of the cell $\Omega$, enclosed by the blue square on the left. On the right, the $\WENO$ interpolator involves the cell vertices and their first neighbours, enclosed in the orange square.}
\label{fig:stencils}	
\end{figure}
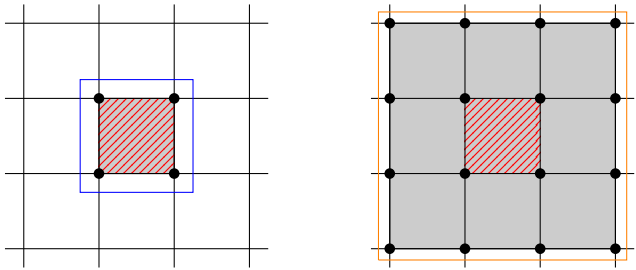

Here instead we resort to two different types of local interpolator with the aim of making the algorithm faster and practical for surface reconstructions on a finer grid, using parallel computing. The first one is a multilinear interpolator, while the second is a $\WENO$ interpolator. The stencils involved in the reconstruction procedures for the 2d case are shown in Fig.~\ref{fig:stencils}.

In this way, we can guarantee that the only communications in the semi-Lagrangian scheme occur when the two (resp. four) interpolation points needed for the update \eqref{eq:SL2d} (respectively \eqref{eq:SL3d}) belong to a different processor than the one owning the point $\vec{x}_j$ (see \S\ref{ssec:parallelization}) and not during the reconstruction procedure itself.

\subsubsection*{Multilinear interpolator}

This interpolator is simply the $\Qone$ Finite Element interpolation on the grid $\grid$.
Let $\{\varphi_j\}_{j\in\grid}$ be the shape functions that are, in each cell, a tensor product of degree 1 polynomials in each space direction and
such that $\varphi_j(\vec{x}_k)=\delta_{j,k}$.
Then, for any function $v(\vec{x})$, we consider the interpolator
\begin{equation} \label{eq:IQ1}
    I_{\Qone}[v](\vec{x}) = \sum_{j\in\grid} v_j \varphi_j(\vec{x})
\end{equation}
where $v_j$ denotes the point value of the function $v$ at the point $\vec{x}_j$.

We remark that, once the 
voxel of the grid containing $\vec{x}$ is located, the summation above in the 2d case (respectively in the 3d case) involves only four (resp. eight) terms and can be locally computed.

We expect that this kind of interpolator will be computationally very cheap, but that it might give rise to artifacts in the reconstructed surfaces since it is piecewise linear and cannot faithfully represent the curvature of the surface.

\subsubsection*{$\WENO$ interpolator}

In order to have a more accurate interpolation operator, we resort to higher order polynomials; the larger stencils involved, however, expose to the risk of introducing oscillations due to the presence of corner points in the data. We thus choose $\WENO$ techniques. With respect to classical $\WENO$ reconstructions employed in conservation laws, here the reconstruction point is arbitrary and not necessarily on a face of the grid.

The procedure for the construction of a high-order $\WENO$ interpolator follows the one described in \cite{CFR05}. For practical purposes, here we will give some details of this procedure focusing on the case of a third-order reconstruction in the cell $\Omega_{j}=[x_j,x_{j+1}]$, on a uniform Cartesian mesh of width $\dx$. Additionally, the description will be limited to the one-dimensional case, since, in the multi-dimensional case, one only needs to iterate this procedure by dimension.

When constructing a $\WENO$ interpolator on the interval $[x_j,x_{j+1}]$, one searches for a convex combination of low-degree polynomials, designed in such a way that a high-order reconstruction is computed in regions associated to smooth data, while non-oscillatory properties are guaranteed in presence of a discontinuity. The selection and blending of such polynomials are delegated to the well-known regularity indicators.

More in details, to construct a $\WENO$ interpolation of degree $2n-1$ on the interval $[x_j,x_{j+1}]$, we start considering the Lagrange polynomial built on the stencil $S=\{x_{j-n+1},\ldots,x_{j+n}\}$, written in the form
\begin{equation}
    Q(x) = \sum_{k=1}^n C_k(x)P_k(x),
\end{equation}
where the $C_k$ are polynomials of degree $n-1$ and the $P_k$ are polynomials of degree $n$ interpolating the point values $v_l$ of a function $v$ on the stencil $S_k=\{x_{j-n+k},\ldots,x_{j+k}\}$, $k=1,\ldots,n$. Then, we proceed as follows:
\begin{enumerate}
    \item compute suitable regularity indicators
    \begin{equation}
        \OSC_k=\OSC[P_k], \quad k=1,\ldots,n;
    \end{equation}
    \item define the quantities
    \begin{equation}
        \alpha_k(x)=\frac{C_k(x)}{(\OSC_k+\epsilon)^2},
    \end{equation}
    with $\epsilon=\dx^2$;
    \item compute the nonlinear weights $\{w_k\}_{k=1}^n$ as
    \begin{equation}
        w_k(x)=\frac{\alpha_k(x)}{\sum_l\alpha_l(x)};
    \end{equation}
    \item finally, define the reconstruction polynomial as
    \begin{equation}\label{IWENO}
        I[V](x) = \sum_{k=1}^n w_k(x)P_k(x).
    \end{equation}
\end{enumerate}

Differently from \cite{CFR05}, we consider regularity indicators defined as
\begin{equation}
\label{eq:ind:1d}
\OSC_k = \OSC[P_k] = \sum_{\beta\geq2} \dx^{2\beta-3} \int_{x_j}^{x_{j+1}} \left(\frac{d^{(\beta)}P_k}{dx^\beta}\right)^2 dx.
\end{equation}
We point out that the above is the classical definition of the oscillation indicators for the $\WENO$ reconstruction as given in \cite{JiangShu:96}, except that the first derivative is not included in the sum. This choice is justified by the fact that the function that we want to interpolate can be at worst continuous with kinks.

Now, turning back to our specific case, we need to find a proper expression for the linear weights $C_k$ in a one-dimensional framework and considering $n=2$ to get, at least in the best case, a third-order reconstruction. Since $n=2$, the whole stencil is given by $S=\{ x_{j-1},x_j,x_{j+1},x_{j+2}\}$ and the only two substencils to consider are $S_L=\{ x_{j-1},x_j,x_{j+1}\}$ and $S_R=\{ x_j,x_{j+1},x_{j+2}\}$. Since the polynomials $P_L$ and $P_R$ interpolate the function $v$ respectively on $S_L$ and $S_R$, each $C_k$ should vanishes outside $S_k$, namely $C_L(x_{j+2})=0$ and $C_R(x_{j-1})=0$, and the condition $\sum_{k=L,R} C_k(x_i)=1$ must hold for every node $x_i\in S$. Thus, we have two conditions for each polynomial $C_k$, which, in our case leads to a unique definition of these one degree polynomials that we can write in the form
\begin{equation} \label{eq:ckgamma}
    C_k(x) = \gamma_k \frac{x-\hat{x}_k}{\dx},
\end{equation}
where $k=L,R$ and $\hat{x}_k \in S\setminus S_k$.

To get the coefficients $\gamma_L$ and $\gamma_R$ we start considering the first node of the stencil $x_{j-1}$, on which we have $C_L(x_{j-1})=1$ and $C_R(x_{j-1})=0$, so that
\begin{equation}
    Q(x_{j-1}) = C_L(x_{j-1})P_L(x_{j-1}) = C_L(x_{j-1}) v_{j-1},
\end{equation}
thus inferring that $C_L(x_{j-1})=1$ and therefore, using \eqref{eq:ckgamma} with $\hat{x}_k=x_{j+2}$, 
\begin{equation}
    \gamma_L = -\frac{1}{3}.
\end{equation}
Analogously, we can proceed considering the node $x_j$ for which we must have $C_L(x_j)+C_R(x_j)=1$. Since the expression for $C_L$ is known, we can compute
\begin{equation}
    C_R(x_j) = 1-C_L(x_j) = 1-\frac{2}{3} = \frac{1}{3}
\end{equation}
and therefore, from $\frac{1}{3}=\gamma_R[(x_j-x_{j-1})/\dx]$ in \eqref{eq:ckgamma}, we get
\begin{equation}
    \gamma_R = \frac{1}{3}.
\end{equation}
Summarizing, we have the following expressions:
\begin{equation}
    C_L(x) = \frac{x_{j+2}-x}{3\dx}, \quad C_R(x) = \frac{x-x_{j-1}}{3\dx}.
\end{equation}

It only remains to find a proper expression of the regularity indicators $\OSC_L$ and $\OSC_R$ for the two degree polynomials $P_L$ and $P_R$. Let us consider a polynomial of degree at most two written in the form $P(x)=\sum_{i=0}^2 a_i[(x-x_j)/\dx]^i$. Its indicator is a quadratic form of its coefficients and is given by
\begin{equation*}
\OSC[P] = \frac{1}{\dx^2}4 a_2^2 
\end{equation*}
or equivalently, denoting by $V$ the vector of data $(v_{j-1},v_j,v_{j+1},v_{j+2})^T$, since the coefficients of the polynomial linearly depends on $V$,
we can express the regularity indicators as a quadratic form of the data being interpolated. For our specific case we obtain
\begin{equation*}
\OSC[P_k] = \frac{1}{\dx^2} V^T A_k V \quad (k=L,R),
\end{equation*}
with
\begin{equation*}
A_L = \begin{pmatrix}
1 & -2 & 1 & 0 \\
-2 & 4 & -2 & 0 \\
1 & -2 & 1 & 0 \\
0 & 0 & 0 & 0 \\
\end{pmatrix} \quad \text{and} \quad
A_R = \begin{pmatrix}
0 & 0 & 0 & 0 \\
0 & 1 & -2 & 1 \\
0 & -2 & 4 & -2 \\
0 & 1 & -2 & 1 \\
\end{pmatrix}.
\end{equation*}

Finally, we illustrate in some details the two-dimensional reconstruction. Let $(x,y)$ be the reconstruction point, located in the cell $[x_i,x_{i+1}]\times[y_j,y_{j+1}]$. One first performs four one-dimensional $\WENO$ interpolations to compute auxiliary data $w_{i-1},w_{i},w_{i+1},w_{i+2}$. Each $w_k$ is the interpolation in the $y$ direction of the data $v_{k,j-1},\ldots,v_{k,j+2}$ evaluated at $y$. Finally a $\WENO$ interpolation in the $x$ direction of $w_{i-1},w_{i},w_{i+1},w_{i+2}$ evaluated at $x$ will be the reconstructed value.
The three-dimensional case can be treated analogously.

Note that, due to the dimensional splitting employed in the two and three dimensional procedures, the $\WENO$ interpolator, as the $\RBF$ and the $\Qone$ ones, is globally continuous on the edges of the cell in which the reconstruction is performed. Moreover, it is more local than $\RBF$ and it is expected to represent curved surfaces more faithfully than $\Qone$. 

\subsection{Domain decomposition for parallel runs}\label{ssec:parallelization}

Especially three-dimensional computations suggest the use of parallel computing for reducing the execution time. Here we recall how the computational grid has been organized for parallel implementation and what are the main communications between workers that might occur during the execution.

As usual for Cartesian grids, in our parallel implementation, each rank $m$ owns a rectangular portion $G_m$ of the computational grid $\grid$ such that $\grid = \cupdot_m G_m$ and is responsible to update the solution on it. Since the computation of the reconstruction and of the time update might require knowledge of the solution at the previous time step in a neighbourhood of $G_m$, each rank also keeps a local copy of the data in a halo region $\widetilde{G}_m$ of at least one ghost point per direction around $G_m$. The width of this halo region depends on which information are needed to perform all the computations in the algorithm. The halo data in $\widetilde{G}_m$ have to be synchronized between workers during the computations through inter-process communications.

In the case of a semi-Lagrangian scheme, the main sources of communication arise from: first, the interpolation of $\phi^n$ at feet of the characteristics emanating backwards from the nodes $\vec{x}_j$ and specified in equations \eqref{eq:SL2d} and \eqref{eq:SL3d}; second, the stencil width required for the computation of the interpolant. 
This latter is mitigated by our choice of local interpolation techniques with very small stencils, like those of \S\ref{ssec:interp}: in both cases we are able to run with one ghost-cell per direction for the $\Qone$ reconstruction and two ghosts per direction in the $\WENO$ case.

The former task is facilitated by the Cartesian structure of the grid and of its partitioning scheme. Consider any of the points $\vec{x}^*_{j,i}$ appearing in \eqref{eq:SL2d} or \eqref{eq:SL3d} associated to a node $\vec{x}_j$ owned by the processor $m$. First $\vec{x}^*_{j,i}$ is located on the grid partition and let $m^*$ be the processor owning the cell into which it falls. If $m^*=m$, the interpolation $I[\phi^n](\vec{x}^*_{j,i})$ is computed without the need of any communication. Otherwise, a request to the processor $m^*$ is sent, communicating the point $\vec{x}^*_{j,i}$ 
and receiving in response the value of $I[\phi^n](\vec{x}^*_{j,i})$. Of course all these communications are gathered in a single push.

In our implementation also the point cloud $\pcloud$ is distributed among the workers. Some more communications are needed in the auxiliary routines described in the next sections and will be highlighted therein.

\subsection{Distance function}\label{ssec:distance}
In the semi-Lagrangian schemes \eqref{eq:SL2d} and \eqref{eq:SL3d} one needs to evaluate the distance $d(\vec{x}) = \min_{Q \in \pcloud}|\vec{x}-Q |$ from the dataset $\pcloud$ at every grid point. 
We remark that a direct computation would have computational complexity 
proportional to $N\times|\pcloud|$, where $N$ is the number of points in the computational grid and $|\pcloud|$ the size of the point cloud. It would also require an overwhelming amount of communications in parallel runs, when both the grid and the point cloud are distributed among different workers.

However, one acknowledges that accurate values of the distance function are needed only close to the object surface to which $\pcloud$ belongs. In practice we initially set $d(\vec{x}_j)$ to the exact distance from $\pcloud$ on all nodes of $\grid$ that are in a box of $4\times4$ or $4\times4\times4$ grid points around each point $Q\in\pcloud$. Then, we apply the fast sweeping method of \cite{Zhao2005AFS} to compute approximate values for $d(\vec{x})$ at all other grid points in $\grid$. These approximate values are still good enough to drive the evolution of the surface $\Gamma$ towards $\pcloud$, while the final stages of the evolution will be guided by the exact values of the distance set in the neighbourhood of the points in $\pcloud$.

\subsection{Energy functional}\label{ssec:energy}
Applying the schemes \eqref{eq:SL2d} and \eqref{eq:SL3d} with $p>1$ also requires to approximate, at each time step, the value of the energy functional $E_p(\phi)$ defined in \eqref{eq:energy:ls}. 
We approximate the Dirac $\delta$ function by restricting the integration domain to the subset 
$\grid_0$ composed by the cells in $\grid$ where, at a specific time step, the front is located. Further, 
we assume that the function $\phi$ is, at least locally around $\grid_0$, a signed distance so that $|\nabla\phi|=1$. This latter point is ensured by the reinitialization procedure described in \S\ref{ssec:reinit}.
We thus compute the energy as
\begin{equation}\label{eq:energyApprox}
    E_p(\phi) \approx \left( \sum_{\Omega\in \grid_0} \int_{\Omega} |d(\vec{x})|^p d\vec{x}\right) ^{1/p}.   
\end{equation}
To detect the cells in $\grid_0$, we check the values of $\phi$ on their vertexes and consider only the cells across which the function $\phi$ changes sign. Once the front is located in a cell $\Omega$, different strategies can be considered to compute \eqref{eq:energyApprox}.
\begin{itemize}
    \item In two space dimensions, we identify 
    two points intercepted by the front of $\phi$ on the cell boundary, by the linear interpolation of $\phi$ on the edges. Then the trapezoidal rule is used to approximate the integral over the segment connecting them. This approach would be much more involved in the three-dimensional case since one would need to identify the (two-dimensional) intersection of the zero level set surface with the cube $\Omega$ and this could be placed in a general position.
    \item In three space dimensions, to avoid the aforementioned complications, we further approximate the Dirac $\delta$ function within the volume $\Omega$ by considering a local refinement of the grid element with $\dx' = \dx/R$ to detect the smaller subcells containing the front, and consider the approximation
    \begin{equation}\label{eq:energyApproxFiner}
        E_p(\phi) \approx \left( \sum_{\Omega'\in \grid'_0} |d(\vec{x}')|^p (\dx')^{n-1} \right) ^{1/p},
    \end{equation}
    where $d(\vec{x}')$ is the interpolated value of the distance function at the center of the subcell $\Omega'$.
    
    The subgrid $\grid'_0$ is composed by the subcells $\Omega'$ such that the reconstruction of $\phi$ at the center of the cell satisfies
    \begin{equation*}
        |\phi(\vec{x}')| < \frac{\sqrt3}{2}\dx'.
    \end{equation*}
    We employ this approach with $R=5$.
\end{itemize}

\subsection{Initial data}\label{subsec:initialData}

The initial data should be chosen as an approximation of the signed distance function from the data set $\pcloud$. Obviously, the better the initial data, the more efficient the method will be, but one has to take into account also the computational effort spent in the computation of the initial datum itself. In practice, we compute $\phi^0$ from the distance function, similarly to \cite{Zhao:2000}, to obtain a signed distance function whose zero level surface encompasses the point cloud and is as close as possible to it. 

In particular we start by the approximate distance function computed as in \S\ref{ssec:distance}.
First, the cells on the outer boundary of $\grid$ are marked as external. Then, moving inwards from all boundaries in all Cartesian directions, we propagate the external point marking to nearby cells until we find a point for which $d(\vec{x}_j)<\gamma_{\pcloud}$, where $\gamma_{\pcloud}$ is a suitable threshold related to the resolution of the point cloud. 

After external regions have been so identified, $\phi^0$ is provisionally set to $d(\vec{x}_j)-\gamma_{\pcloud}$ on external points and to an artificially high value otherwise. Next, the fast sweeping method of \cite{Zhao2005AFS} is applied to $\phi^0$, recomputing its values at internal points; finally, the sign of $\phi^0$ is reversed on internal points.

\begin{figure}
\begin{minipage}{0.33\linewidth}
    \centering
    \includegraphics[width=0.8\linewidth]{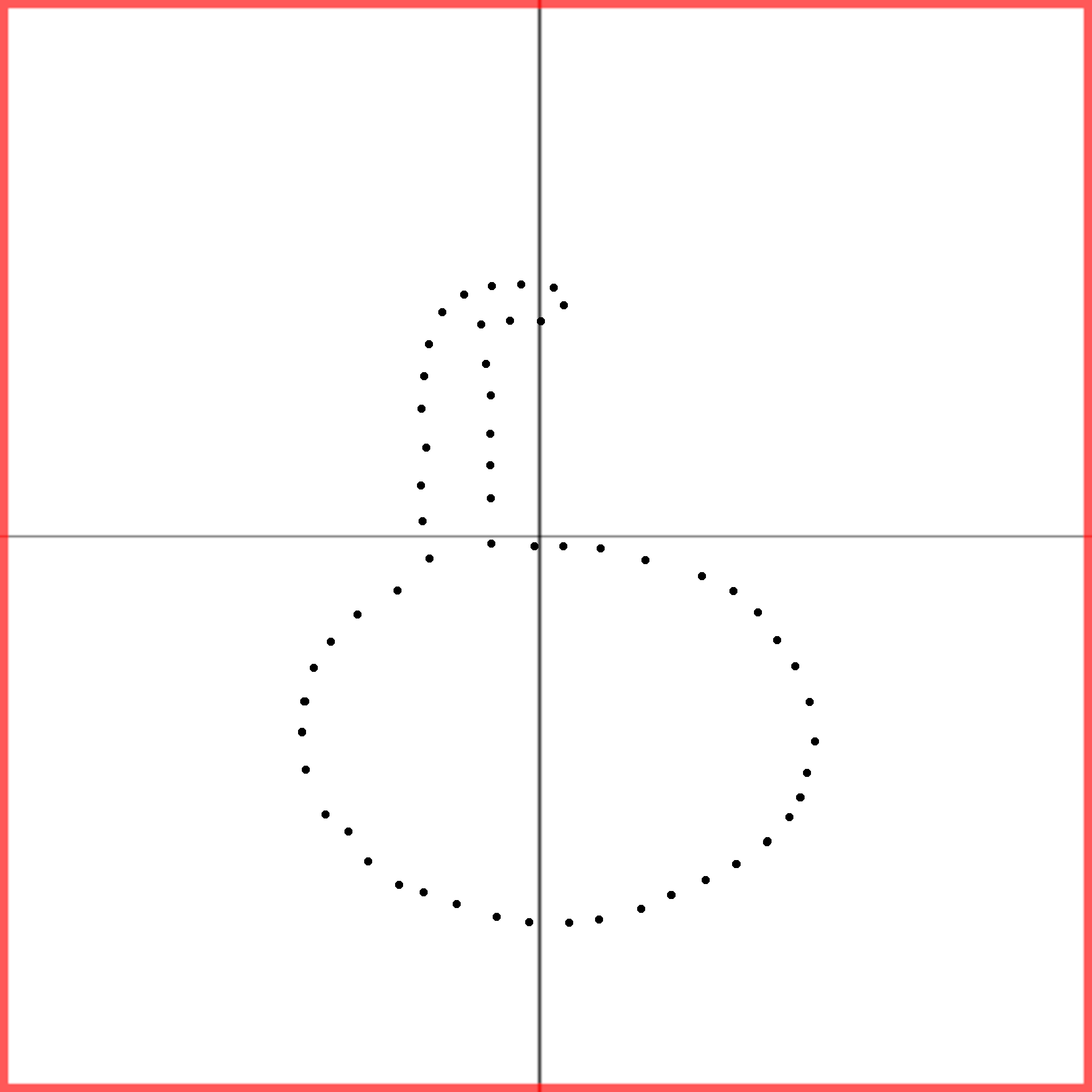}
\end{minipage}\hfill
\begin{minipage}{0.33\linewidth}
    \centering
    \includegraphics[width=0.8\linewidth]{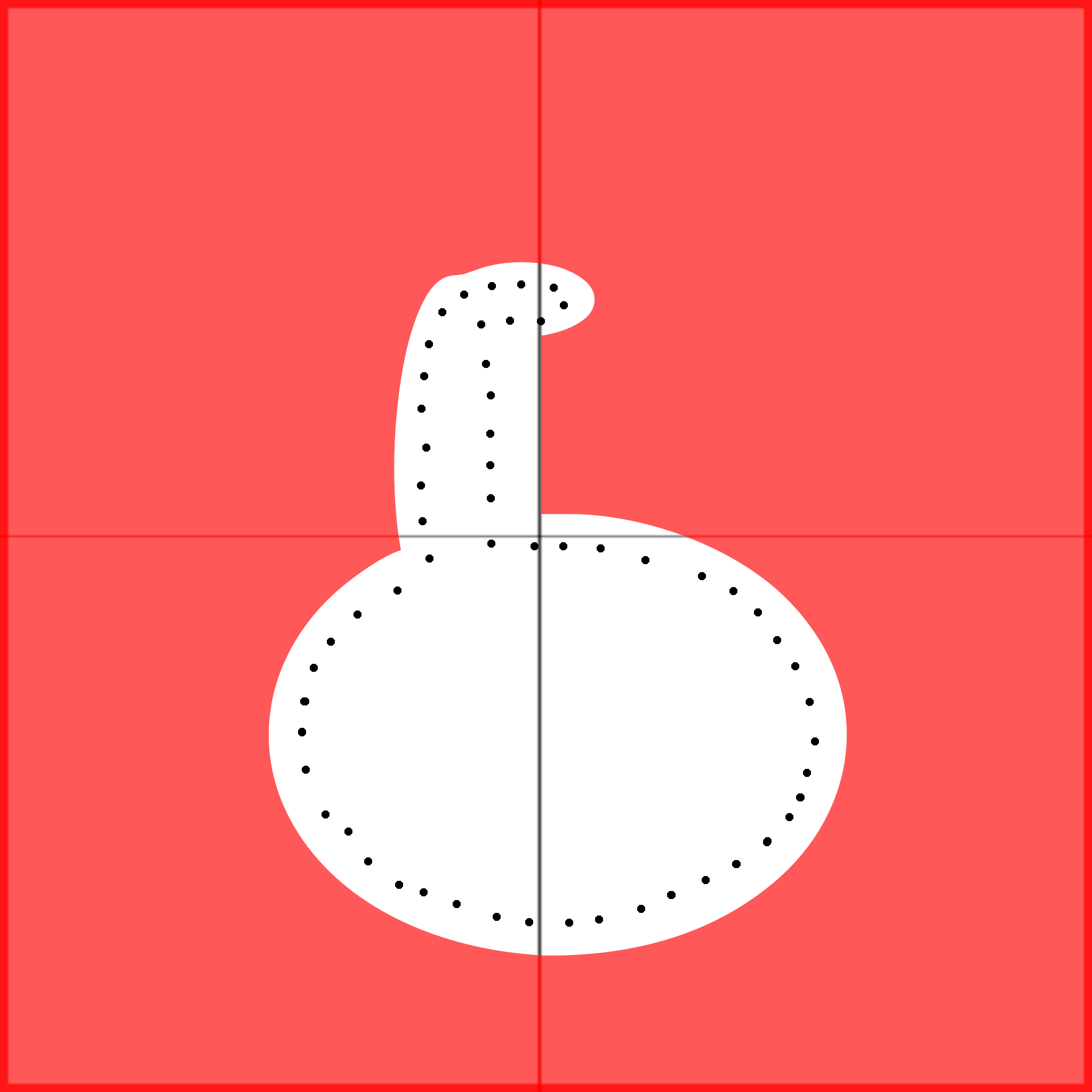}
\end{minipage}\hfill
\begin{minipage}{0.33\linewidth}
    \centering
    \includegraphics[width=0.8\linewidth]{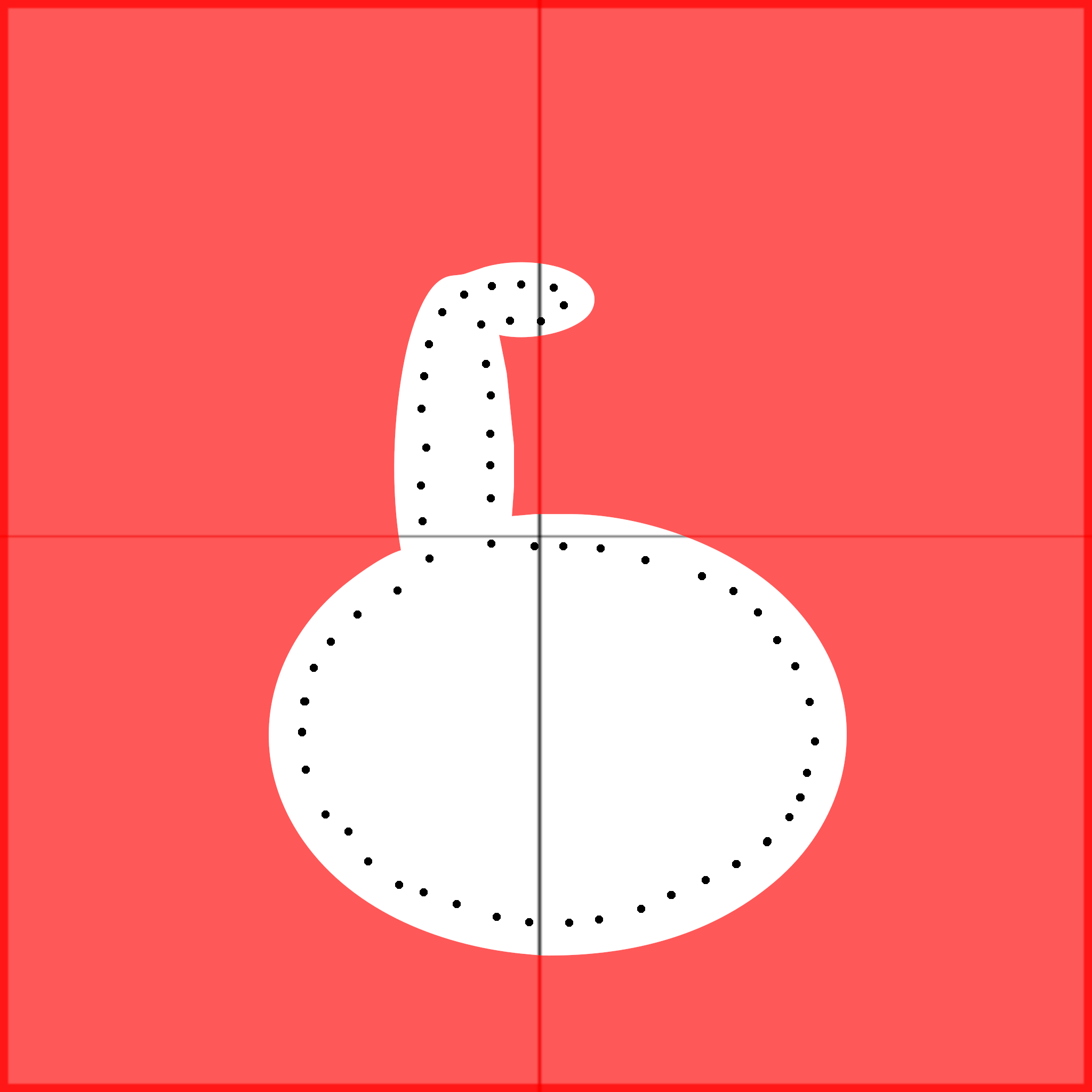}
\end{minipage}
    \caption{Illustration of marking external points in parallel run with a $2\times2$ domain decomposition (grey lines). Left: initial state. Center: after the first sweep. Right: after the second sweep.}
    \label{fig:externalMarking}
\end{figure}

We point out that in parallel runs the procedure for initial marking of external points is similar to the serial run, but is iterated more times and interleaved with communications of point marking in the halo regions.
This is illustrated in Fig.~\ref{fig:externalMarking}. The top-left processor, in the first run (middle panel), expands the marking in the right and bottom direction from the outer border of the computational grid and is thus unable to mark
the small area under the hook. In the second sweep (right panel), the external marking is propagated from the marking done by the top-right processor.
A number of iterations of at least half the number of processors per direction in the decomposition ensures that the initial marking obtained is the same that would have been computed in a serial run.

The one described above, is a good strategy to find an initial data and start the evolution having only the information carried by the point cloud. As a matter of fact, one could also choose to evolve a less accurate initial guess until steady-state on a coarse grid and then use this final level set function as the initial data for a finer evolution.  

\subsection{Reinitialization}\label{ssec:reinit}
One of the nice feature of using a signed distance function to capture interfaces is that some geometric quantities, such as normals and curvatures, are easier to compute in terms of $\phi$. In particular this is of paramount importance when the level set is used in ghost-fluid algorithms for the discretization of PDEs. The evolution determined by \eqref{eq:levelset:pde} will push the solution towards the data set, but will create sharp gradients in the level set function $\phi$ around it.
To keep $|\nabla\phi|$ close to $1$, we employ a reinitialization procedure, first introduced in \cite{SSO:1994}, which consists in evolving until steady-state the equation
\begin{equation}\label{eq:reinit}
\begin{cases}
\pder{\phi}{\tau} + \mathrm{sign}(\tilde\phi) \left( |\nabla\phi| -1\right) =0
\\
\phi(\vec{x},0) = \tilde\phi(\vec{x})
\end{cases}
\end{equation}
where the initial value $\tilde\phi$ corresponds to the solution of \eqref{eq:levelset:pde} at a fixed time step and $\tau$ is a pseudo-time.

The delicate issue with the reinitialization is to guarantee that the interface of $\tilde\phi$ is displaced as little as possible, ideally not at all, while the level set function is modified to ensure that $|\nabla \phi| = 1$. To this aim, we have followed the constrained reinitialization (CR) scheme introduced by Hartmann et al. \cite{HaMeSc:2008} which consists in a modification of the first scheme by Sussman et al. \cite{OshSet:1988} in order to locally preserve the location of the interface during the reinitialization procedure. 

With this approach a one step procedure is performed to reinitialize the level set function in the first neighbours of the interface, which can be easily detected from the change of sign of $\tilde\phi$.
After that, an iterative procedure is performed to reinitialize the values of $\tilde\phi$ in a proper tube around the already reinitialized region; to this end a first-order spatial discretization and forward Euler integration in pseudo-time $\tau$ is used to approximate the solution of the original equation \eqref{eq:reinit}.
In parallel runs, each step of this algorithm also requires an update of the $\phi$ data in the halo region.

\subsection{Localizing the computational effort}\label{ssec:mask}
A typical approach when dealing with level set methods is to localize the evolution only in a narrow band around the zero level set of $\phi$. Since it is not important to update $\phi$ far away from that, at each time step we choose a subgrid
\begin{equation}\label{eq:subgrid}
\widetilde\grid=\{\vec{x}_j\in\grid: \phi^n_j<\gamma\}\subset\grid,     
\end{equation}
with $\gamma=4\dx$ for the 2d case and $\gamma=6\dx$ for the 3d case, as suggested in \cite{PENG1999}, and update $\phi^{n+1}$ only therein. Outside $\widetilde\grid$ we simply cut our level set function as
\begin{equation}\label{eq:cut}
    \phi(\Vec{x}) = 
    \begin{cases}
        \gamma \quad &\text{if} \quad \phi(\vec{x})>\gamma,\\
        \phi(\vec{x}) \quad &\text{if} \quad |\phi(\vec{x})|\leq\gamma,\\
        -\gamma \quad &\text{if} \quad \phi(\vec{x})<-\gamma.
    \end{cases} 
\end{equation}

To prevent numerical oscillations at the boundary of $\widetilde\grid$, we also update the solution involving an additional cut-off function $c(\phi)$, the same described in \cite{PENG1999}, which is given by
\begin{equation}\label{eq:cutoff}
    c(\phi)=
    \begin{cases}
        1 \hfill &\text{if} \; |\phi|\leq\beta,\\
        (|\phi|-\gamma)^2(2|\phi|+\gamma-3\beta)/(\gamma-\beta)^3 &\text{if} \; \beta<|\phi|\leq\gamma,\\
        0 \hfill &\text{if} \; |\phi|>\gamma,
    \end{cases}
\end{equation}
with $\beta=2\dx$ for the 2d case and $\beta=3\dx$ for the 3d case, thus considering the modified equation
\small
\begin{equation}\label{eq:ZhaoPDECutOff}\nonumber
    \pder{\phi}{t}(\vec{x},t) = c(\phi) \Bigg[  \frac{d(\vec{x})}{E_p(\phi)}   \Bigg]^{p-1} \Bigg( \nabla d(\vec{x})\cdot \nabla \phi(\vec{x},t) + \frac{d(\vec{x})}{p} \nabla\cdot\left( \frac{\nabla \phi(\vec{x},t)}{|\nabla \phi(\vec{x},t)|} \right)|\nabla \phi(\vec{x},t)| 
      \Bigg)
\end{equation}
\normalsize
and properly modifying the semi-Lagrangian schemes \eqref{eq:SL2d} and \eqref{eq:SL3d} to incorporate the additional factor $c(\phi)$. 

\begin{figure}
    \begin{minipage}{0.45\textwidth}
        \centering
        \includegraphics[width=0.9\textwidth]{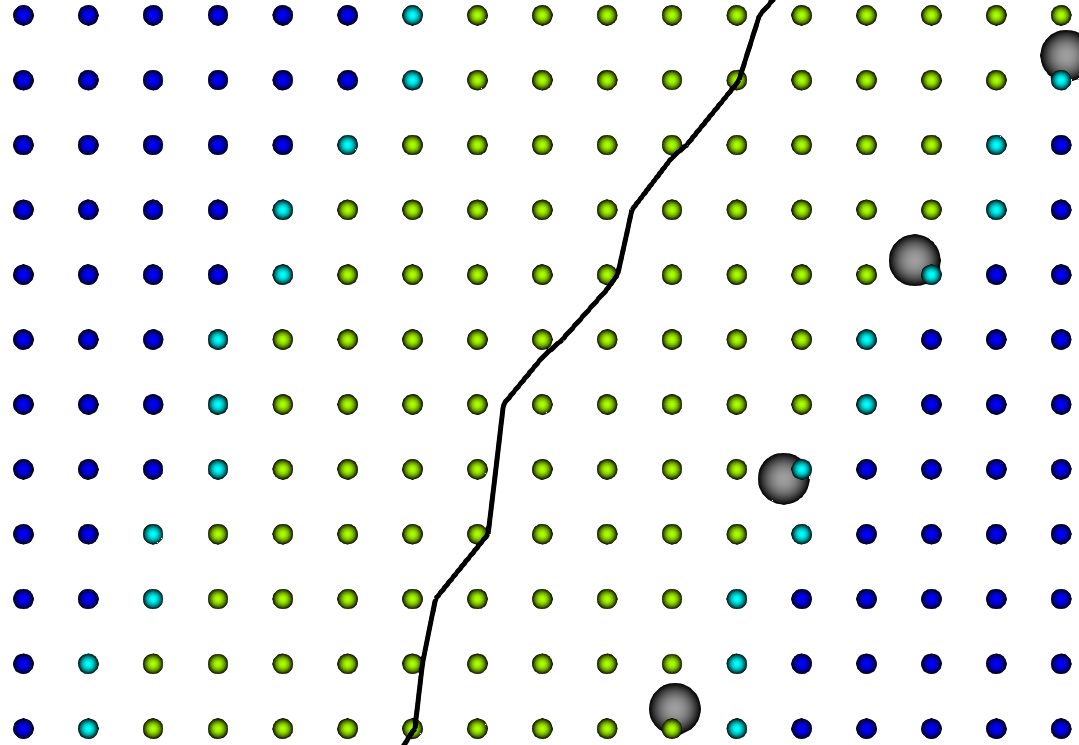}    
    \end{minipage}
    \hfill
    \begin{minipage}{0.45\textwidth}
        \centering
        \includegraphics[width=0.9\textwidth]{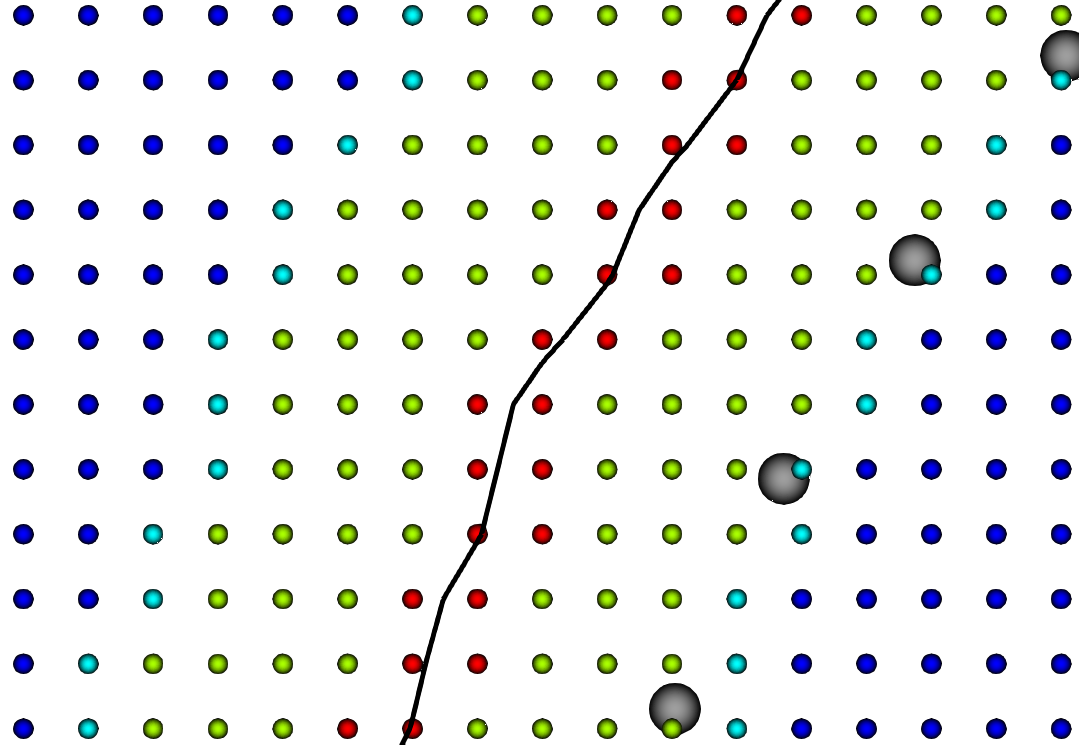}    
    \end{minipage}
    \caption{On the left, the mask computed to detect the computational subgrid $\widetilde\grid$.
    Active nodes are depicted in green, their first neighbours are depicted in light blue (they are inactive during the update, while they are active during the reinitialization step) and remaining inactive nodes are depicted in blue. On the right, the mask computed to detect the computational subgrid $\overline{\grid}$ 
    of the reinitialization. Colours are used in the same way, with in addition red nodes representing the nodes immediately close to the interface on which the signed distance function is computed explicitly with the one step procedure.} 
    \label{fig:mask}
\end{figure}

\begin{figure}
    \centering
    \includegraphics[width=0.6\textwidth]{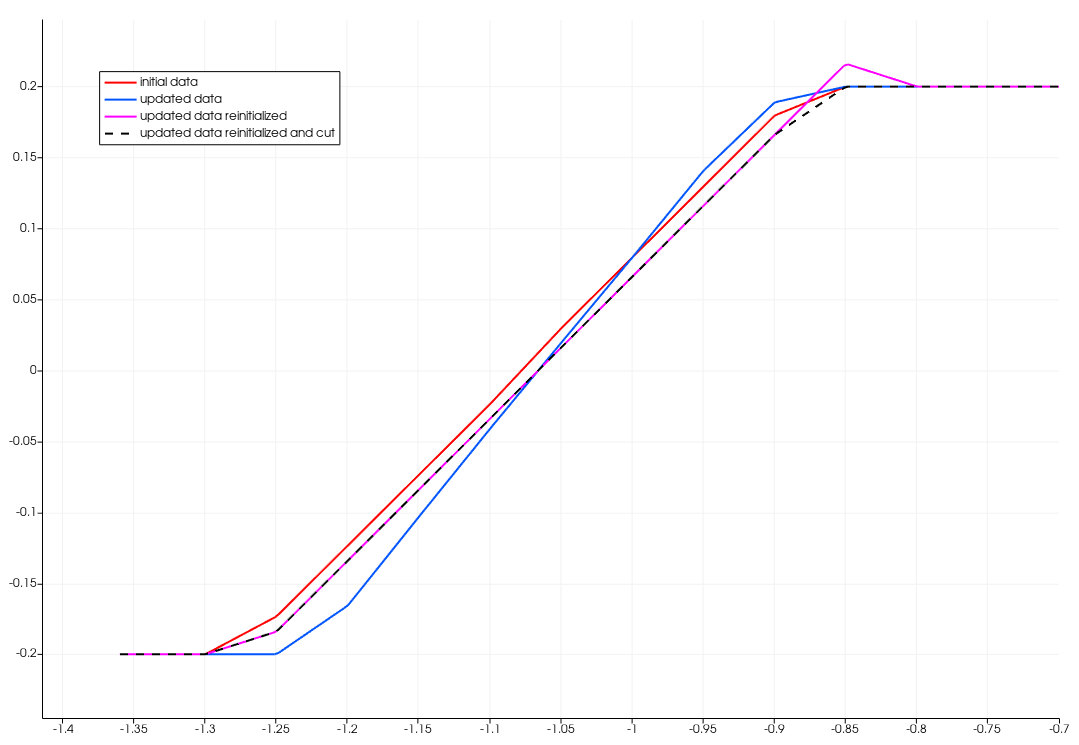}
    \caption{Typical evolution of a front in 1d. Starting from the red line, the update makes the front move towards the data resulting in the blue line, sharper than the red one. Reinitialization (pink line) makes the computational narrow band move contextually with the front, fixing the zero level set. The dotted black line represents the last cutting step necessary to obtain the final update.}
    \label{fig:front}
\end{figure}

Finally, the reinitialization procedure is performed on a different narrow band $\overline\grid$ obtained by considering $\widetilde\grid$ and its first neighbours. The fact that  $\overline\grid \supset \widetilde\grid$ is essential if one starts the algorithm with an initial data so far from the point cloud that the initial computational band \eqref{eq:subgrid} does not contain $\pcloud$. In such a case, without the enlarged reinitialization band, the evolution of $\phi$ would remain confined to the first computational band $\widetilde\grid^0$, while with this choice, the successive bands $\widetilde\grid^n$ will be able to move contextually with the zero level set of $\phi^n$. For more clarity, the bands involved in the different parts of the scheme are depicted in a 2d case in Fig.~\ref{fig:mask} and a typical evolution of a front of the level set function in 1d is depicted in Fig.~\ref{fig:front}.

\section{Numerical experiments}
\label{sec:numtests}

In this section we present some representative numerical results obtained with the algorithm described in the previous sections. As we are looking for a minimum of the energy functional \eqref{eq:energy}, we will consider the same stopping criterion used in \cite{He:2019}: at time $n$, the algorithm stops if
\begin{equation}\label{eq:deltaEn}%
    \Delta_E^n = \frac{|\Bar{e}^k_{n-1}-\Bar{e}^k_{n}|}{\Bar{e}^k_{n}} < 10^{-4}, \quad \text{where} \quad \Bar{e}^k_{n} = \frac{1}{k}\sum_{i=n-k+1}^n E_2(\phi^i)
\end{equation}
or after a maximum of $100$ iterations. Note that in \eqref{eq:deltaEn}, the energy functional \eqref{eq:energy} is computed with $p=2$. This condition is in practice a way to detect stationary points or flat areas of the energy functional. In all the tests we set $k=min(n,10)$ and forced the algorithm to do at least $10$ iterations.

In order to better analize the performance of the algorithm, we also compute, alongside \eqref{eq:deltaEn}, also the $L^1$-norm of the update between two successive iterations and the $L^1$-norm of the error, when the exact level set is given.
Furthermore, we compute the average of the error on the points of the cloud  
\begin{equation}\label{errCloud}
    Err_{\pcloud}^n = \frac{\sum_{\vec{x}_j\in \pcloud} | I[\phi^n]\left( \vec{x}_j\right)|} {|\pcloud|}
\end{equation}
to make some considerations on the role of the curvature regularization and to evaluate how much the final reconstruction is attached to the data set.
\newline

It only remains now to detail the choices for the parameters $p$ and $\delta$ and for the spatial and temporal discretization. 

In \cite{CaFe:2017}, the authors consider the case with $p=1$, which grants for a faster evolution of the initial data towards the cloud and also does not require to compute the factor $[d(x)/E_p]^{p-1}$, while in \cite{Zhao:2000} the authors suggest to set $p=2$, which is more effective in reaching a steady state for the evolution. Also, setting $\delta=0$ simplifies the update formulas \eqref{eq:SL2d} and \eqref{eq:SL3d} since all $\vec{x}^*_{j,i}$ coincide. Furthermore, disregarding the curvature effect, it prevents from loosing too much details, especially on coarse grids. On the other hand, increasing values of $\delta\leq 1$ will smoothen the solution, at the price of having a zero level set slightly off from the cloud $\pcloud$. 
Our choice will thus consist in combining these different approaches in more than one run of Algorithm~\ref{algo}, gradually increasing the resolution of the grid.
\begin{table}[]
    \centering
    \begin{tabular}{|P{0.8cm}|P{0.8cm}|P{0.8cm}|}
    \hline
        $r$ & $p$ & $\delta$ \\\hline
        $1$ & $1$ & $0$ \\
        $2$ & $2$ & $0$ \\
        $\geq3$ & $2$ & $1$ \\\hline
    \end{tabular}
    \caption{Setting of parameters $p$ and $\delta$ in each run $r$ of Algorithm~\ref{algo}.}
    \label{tab:param}
\end{table}
Let $r$ represent the number of the run of Algorithm~\ref{algo}, Table~\ref{tab:param} summarizes our choices for the parameters $p$ and $\delta$, while for spatial and temporal discretization we set
\begin{equation}\label{eq:dxdt}
    \dx^{(r)} = \frac{1}{2^{r-1}}h_{\mathcal{S}}, \quad \dt^{(r)} = \frac{1}{4}\dx^{(r)}.
\end{equation}
Unless otherwise specified, three runs of Algorithm~\ref{algo} are performed.

In \eqref{eq:dxdt} $h_{\mathcal{S}}$ represents an estimate of the resolution of the cloud $\pcloud$; its value is approximated during preprocessing by randomly choosing a sample made up of the $10\%$ of the points in $\pcloud$ and then computing the average of the distances between each of these points and their nearest neighbour in $\pcloud$.

In the first and coarsest run ($r=1$), we compute the initial data as described in subsection \ref{subsec:initialData}. The required threshold is set as
$\gamma_{\mathcal{S}} = K_{\mathcal{S}}\, h_{\mathcal{S}}$
with $K_{\mathcal{S}}$ set to $2.0$ unless otherwise specified. 
The role of $K_{\mathcal{S}}$ is relevant when the points in $\pcloud$ are not evenly distributed on the sampled surface. It is in fact common, especially in 3d, due to the  supports used for the object during laser scanning, to find data sets that have piggy bank-like shapes, which have fake holes that could distort the result. In such cases, the good practice consists in choosing a larger $K_{\mathcal{S}}$, to start with level set that is further away from the data, but encloses the point cloud without entering in the fake cavities. For the later runs ($r>1$), we use as initial data the steady-state solution computed by the previous run and interpolated on the current grid.%

All the codes has been written in C++ language with the support of the MPI and of the PETSc library \cite{petsc-user-ref,petsc-efficient}. In particular, the communications in the semi-Lagrangian step \ref{algo:sl} of Algorithm~\ref{algo} has been performed with the help of a DMSWARM object.
The tests has been performed on the Galileo100 cluster hosted at CINECA\footnote{https://www.hpc.cineca.it/systems/hardware/galileo100/}.

\subsection{2d data sets}

\subsubsection{Circle test}
We first consider a data set $\pcloud$ made up of $64$ points uniformly chosen on a circle of radius $1$, thus having a cloud size approximately equal to $9.81\times10^{-2}$. The preprocessing provides the expected approximation of $h_{\pcloud}$ and sets $\gamma_{\pcloud} = 9.81\times10^{-2}$.

\begin{figure}
\begin{minipage}{0.33\linewidth}
    \centering
    \includegraphics[width=1.0\linewidth]{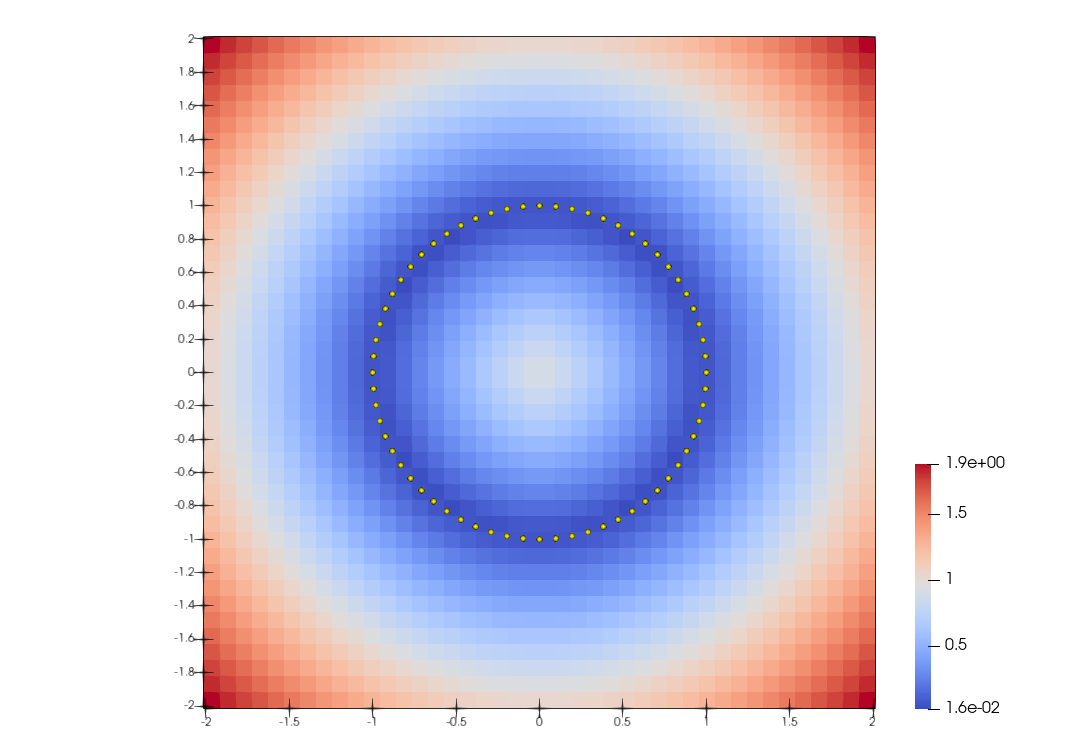}
\end{minipage}\hfill
\begin{minipage}{0.33\linewidth}
    \centering
    \includegraphics[width=1.0\linewidth]{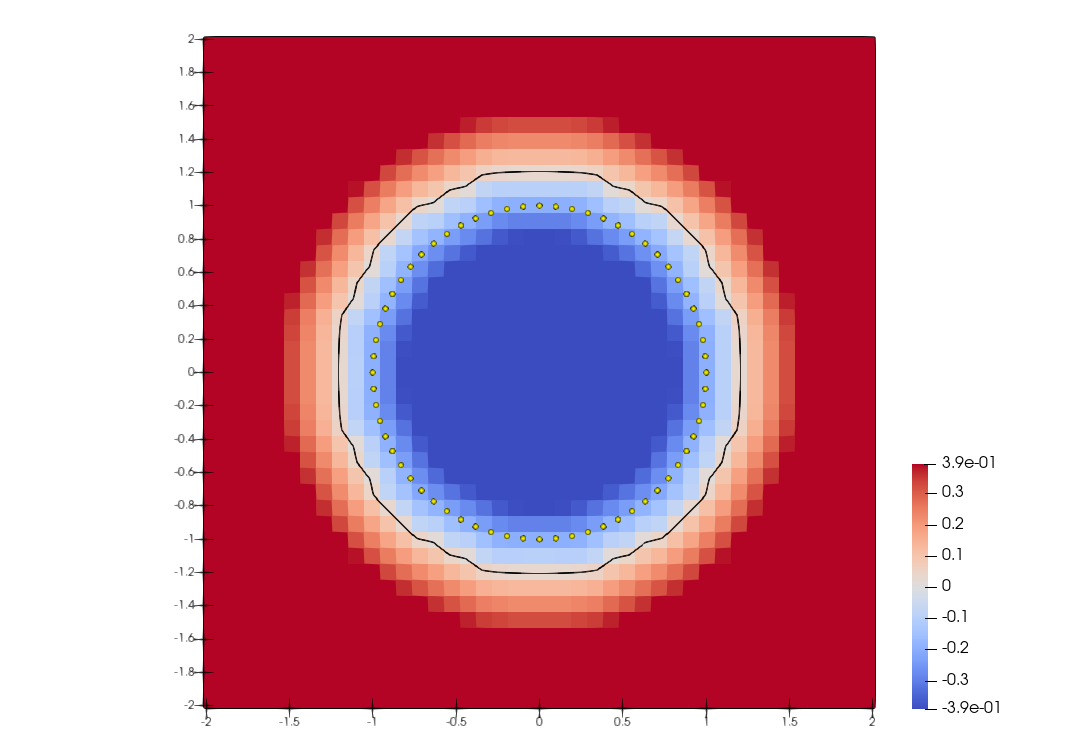}
\end{minipage}\hfill
\begin{minipage}{0.33\linewidth}
    \centering
    \includegraphics[width=1.0\linewidth]{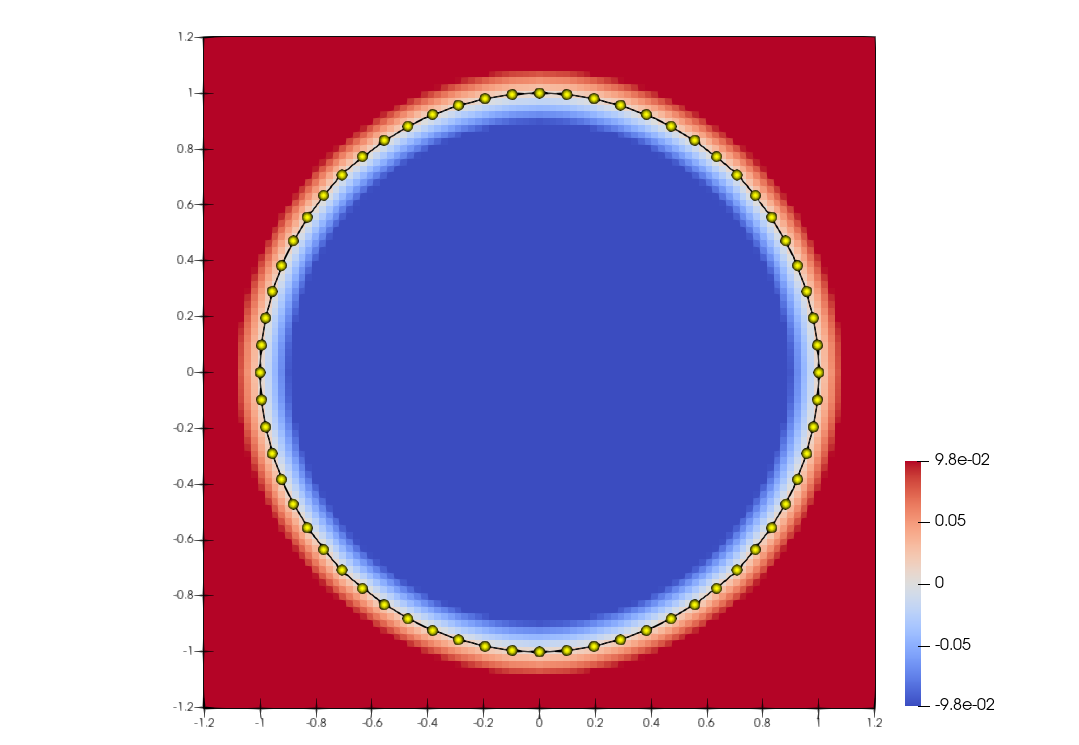}
\end{minipage}
    \caption{Steps of the algorithm for the 2d circle case: on the left, the distance function is represented; the central panel shows the initial data and its zero level set (black line); the final data and its contour (black line) are represented on the right. Data represented here has been obtained with $\WENO$ reconstruction.}
    \label{fig:circleReconstruction}
\end{figure}

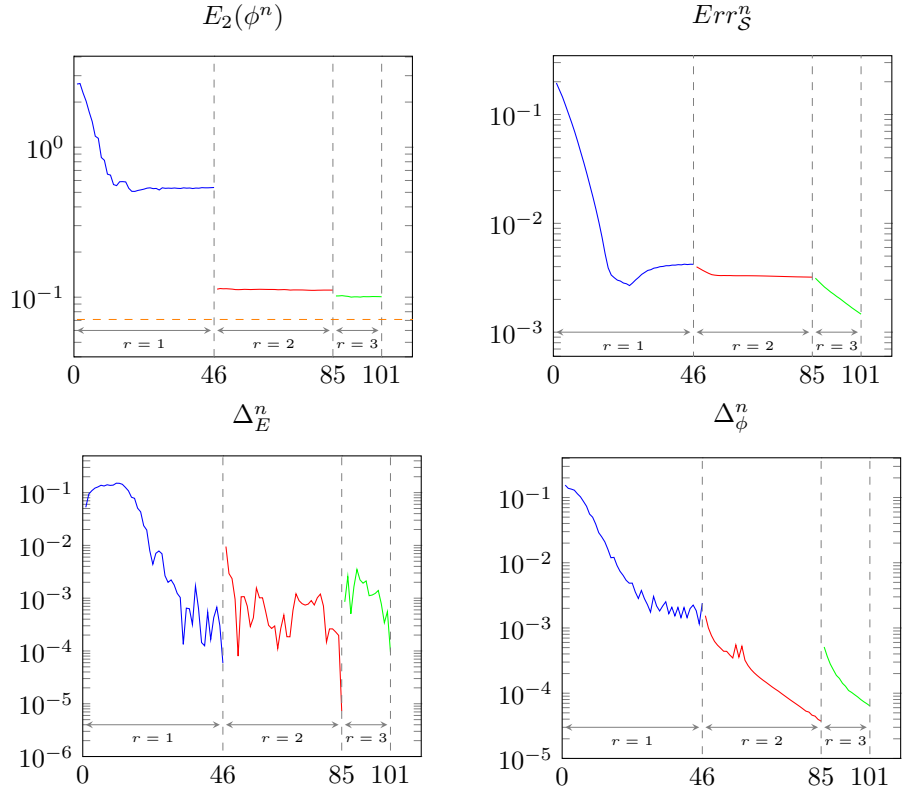
\begin{figure}
\begin{minipage}{0.45\linewidth}
    \begin{tikzpicture}

\begin{semilogyaxis}[%
  width=4.5cm,height=4cm,scale only axis,
  xtick={0, 46, 85, 101},
  xticklabels={$0$, $46$, $85$, $101$},
  major tick length = 0,
  yticklabel style={/pgf/number format/sci},
  ymin=4e-2,xmin=0,
  legend columns=-1, 
  legend style={
  	legend pos=north east,
  	draw=none,
  },
  title={$E_2(\phi^n)$}
]
\addplot[color=blue] table[col sep=space]{grafici/circleEn1.txt};
\addplot[color=red] table[col sep=space]{grafici/circleEn2.txt};
\addplot[color=green] table[col sep=space]{grafici/circleEn3.txt};

\draw[dashed,help lines] (axis cs:46,1e-20) -- (axis cs:46,100) ;
\draw[dashed,help lines] (axis cs:85,1e-20) -- (axis cs:85,100) ;
\draw[dashed,help lines] (axis cs:101,1e-20) -- (axis cs:101,100) ;
\draw[stealth-stealth,help lines] (axis cs:1,6e-2) -- node[pos=0.5,below,black]{\tiny $r=1$} (axis cs:45,6e-2) ;
\draw[stealth-stealth,help lines] (axis cs:47,6e-2) -- node[pos=0.5,below,black]{\tiny $r=2$} (axis cs:84,6e-2) ;
\draw[stealth-stealth,help lines] (axis cs:86,6e-2) -- node[pos=0.5,below,black]{\tiny $r=3$} (axis cs:100,6e-2) ;

\draw[dashed,help lines,orange] (axis cs:1,7.103510780002126e-02) -- (axis cs:150,7.103510780002126e-02) ;

\end{semilogyaxis}
\end{tikzpicture}
\end{minipage}\hspace{0.8cm}
\begin{minipage}{0.45\linewidth}
    \begin{tikzpicture}

\begin{semilogyaxis}[%
  width=4.5cm,height=4cm,scale only axis,
  xtick={0, 46, 85, 101},
  xticklabels={$0$, $46$, $85$, $101$},
  major tick length = 0,
  yticklabel style={/pgf/number format/sci},
  ymin=6e-4,
  xmin=0,
  legend columns=-1, 
  legend style={
  	legend pos=north east,
  	draw=none,
  },
  title={$Err_{\pcloud}^n$}
]
\addplot[color=blue] table[col sep=space]{grafici/circleErrCloud1.txt};
\addplot[color=red] table[col sep=space]{grafici/circleErrCloud2.txt};
\addplot[color=green] table[col sep=space]{grafici/circleErrCloud3.txt};

\draw[dashed,help lines] (axis cs:46,1e-20) -- (axis cs:46,100) ;
\draw[dashed,help lines] (axis cs:85,1e-20) -- (axis cs:85,100) ;
\draw[dashed,help lines] (axis cs:101,1e-20) -- (axis cs:101,100) ;
\draw[stealth-stealth,help lines] (axis cs:1,1e-3) -- node[pos=0.5,below,black]{\tiny $r=1$} (axis cs:45,1e-3) ;
\draw[stealth-stealth,help lines] (axis cs:47,1e-3) -- node[pos=0.5,below,black]{\tiny $r=2$} (axis cs:84,1e-3) ;
\draw[stealth-stealth,help lines] (axis cs:86,1e-3) -- node[pos=0.5,below,black]{\tiny $r=3$} (axis cs:100,1e-3) ;

\end{semilogyaxis}
\end{tikzpicture}
\end{minipage}\\
\vspace{0.8cm}
\noindent
\begin{minipage}{0.45\linewidth}
    \begin{tikzpicture}

\begin{semilogyaxis}[%
  width=4.5cm,height=4cm,scale only axis,
  xtick={0, 46, 85, 101},
  xticklabels={$0$, $46$, $85$, $101$},
  major tick length = 0,
  yticklabel style={/pgf/number format/sci},
  ymin=1e-6,
  xmin=0,
  legend columns=-1, 
  legend style={
  	legend pos=north east,
  	draw=none,
  },
  title={$\Delta_E^n$}
]
\addplot[color=blue] table[col sep=space]{grafici/circleMeanEn1.txt};
\addplot[color=red] table[col sep=space]{grafici/circleMeanEn2.txt};
\addplot[color=green] table[col sep=space]{grafici/circleMeanEn3.txt};

\draw[dashed,help lines] (axis cs:46,1e-20) -- (axis cs:46,100) ;
\draw[dashed,help lines] (axis cs:85,1e-20) -- (axis cs:85,100) ;
\draw[dashed,help lines] (axis cs:101,1e-20) -- (axis cs:101,100) ;
\draw[stealth-stealth,help lines] (axis cs:1,4e-6) -- node[pos=0.5,below,black]{\tiny $r=1$} (axis cs:45,4e-6) ;
\draw[stealth-stealth,help lines] (axis cs:47,4e-6) -- node[pos=0.5,below,black]{\tiny $r=2$} (axis cs:84,4e-6) ;
\draw[stealth-stealth,help lines] (axis cs:86,4e-6) -- node[pos=0.5,below,black]{\tiny $r=3$} (axis cs:100,4e-6) ;

\end{semilogyaxis}
\end{tikzpicture}
\end{minipage}\hspace{0.8cm}
\begin{minipage}{0.45\linewidth}
    \begin{tikzpicture}

\begin{semilogyaxis}[%
  width=4.5cm,height=4cm,scale only axis,
  xtick={0, 46, 85, 101},
  xticklabels={$0$, $46$, $85$, $101$},
  major tick length = 0,
  yticklabel style={/pgf/number format/sci},
  ymin=1e-5,
  xmin=0,
  legend columns=-1, 
  legend style={
  	legend pos=north east,
  	draw=none,
  },
  title={$\Delta_\phi^n$}
]
\addplot[color=blue] table[col sep=space]{grafici/circleVar1.txt};
\addplot[color=red] table[col sep=space]{grafici/circleVar2.txt};
\addplot[color=green] table[col sep=space]{grafici/circleVar3.txt};

\draw[dashed,help lines] (axis cs:46,1e-20) -- (axis cs:46,100) ;
\draw[dashed,help lines] (axis cs:85,1e-20) -- (axis cs:85,100) ;
\draw[dashed,help lines] (axis cs:101,1e-20) -- (axis cs:101,100) ;
\draw[stealth-stealth,help lines] (axis cs:1,3e-5) -- node[pos=0.5,below,black]{\tiny $r=1$} (axis cs:45,3e-5) ;
\draw[stealth-stealth,help lines] (axis cs:47,3e-5) -- node[pos=0.5,below,black]{\tiny $r=2$} (axis cs:84,3e-5) ;
\draw[stealth-stealth,help lines] (axis cs:86,3e-5) -- node[pos=0.5,below,black]{\tiny $r=3$} (axis cs:100,3e-5) ;

\end{semilogyaxis}
\end{tikzpicture}
\end{minipage}\\

    \caption{Significant quantities computed for the circle test using $\WENO$ interpolator. In the top-left panel the energy functional is depicted: the drop between the first and the second run is due to the change of the grid size and the consequent new computation of the distance function. In the bottom-right panel the error computed on the cloud is shown: as expected, its profile is not monotone since during the evolution the interface can even pass the cloud and then come back. The bottom-left panel represents the running average of the energy functional: note that in the first $10$ iterations of each run the moving average is still forming up. The bottom-right panel represents the $L^1$-norm of the variation of $\phi$ between two successive runs.}
    \label{fig:error_cerchio2}
\end{figure}

\begin{table}
\begin{center}
\tiny
\begin{tabular}{|c|c|c|c|c|c|c|c|}
\hline & & \multicolumn{3}{|c|}{$\Qone$} &\multicolumn{3}{|c|}{$\WENO$} \\
\hline $r$ & Grid & $L^{1}$-err & Energy & Error on $\pcloud$ & $L^{1}$-err & Energy & Error on $\pcloud$ \\ \hline \hline
$1$ & $42\times42$ & $3.51e-02$ & $1.83e-01$ & $4.17e-03$ & $3.81e-02$ & $1.59e-01$ & $4.20e-03$\\ \hline 
$2$ & $58\times58$ & $8.34e-03$ & $1.16e-01$ & $3.03e-03$ & $8.26e-03$ & $1.12e-01$ & $3.20e-03$\\ \hline 
$3$ & $99\times99$ & $2.27e-03$ & $1.03e-01$ & $1.28e-03$ & $2.08e-03$ & $1.00e-01$ & $1.46e-03$\\ \hline 
$4$ & $181\times181$ & $9.45e-04$ & $1.04e-01$ & $1.01e-03$ & $8.28e-04$ & $1.02e-01$ & $9.76e-04$\\ \hline 
$5$ & $344\times344$ & $4.21e-04$ & $1.04e-01$ & $9.26e-04$ & $3.46e-04$ & $1.03e-01$ & $7.90e-04$\\ \hline 
\end{tabular}
\end{center}
\caption{Errors and energy computed for the circle test at the end of each run. To better appreciate the results Algorithm~\ref{algo} has been performed five times, fixing the parameters according to Table~\ref{tab:param}.}
\label{tab:errGlobCircle}
\end{table}

The main steps of the reconstruction of the circle are depicted in Fig~\ref{fig:circleReconstruction}, while Fig~\ref{fig:error_cerchio2} collects significant graphs to evaluate the accuracy of the method. The three runs of the algorithm required respectively $27$, $17$ and $21$ iterations using the $\Qone$ reconstruction, while $46$, $39$ and $16$ iterations are required in the $\WENO$ case, with an amount of computational time equal to $1.78e-01$ and $5.44e-01$ seconds, respectively. Looking at Fig.~\ref{fig:error_cerchio2} note how the zero level set approaches the cloud: as expected, the $L^1$-norm of the variation do not show a monotone profile and the error on the cloud $Err^n_{\pcloud}$ is not vanishing due to the finite grid size and to the curvature term.
During the first run, the interface quickly moves towards the data, it can even pass through the cloud, and then stabilizes. During the successive runs the reconstruction becomes more accurate: the energy  gets closer to the exact value of the energy computed for a circle, represented by the dashed orange line ($\approx 7.10\times10^{-2}$), and the error made on the cloud decreases. In this case, the exact signed distance function from the circle is known and we can compute the global error (see Table~\ref{tab:errGlobCircle}).

\subsubsection{Square test}

\begin{figure}
\begin{minipage}{0.33\linewidth}
    \centering
    \includegraphics[width=1.0\linewidth]{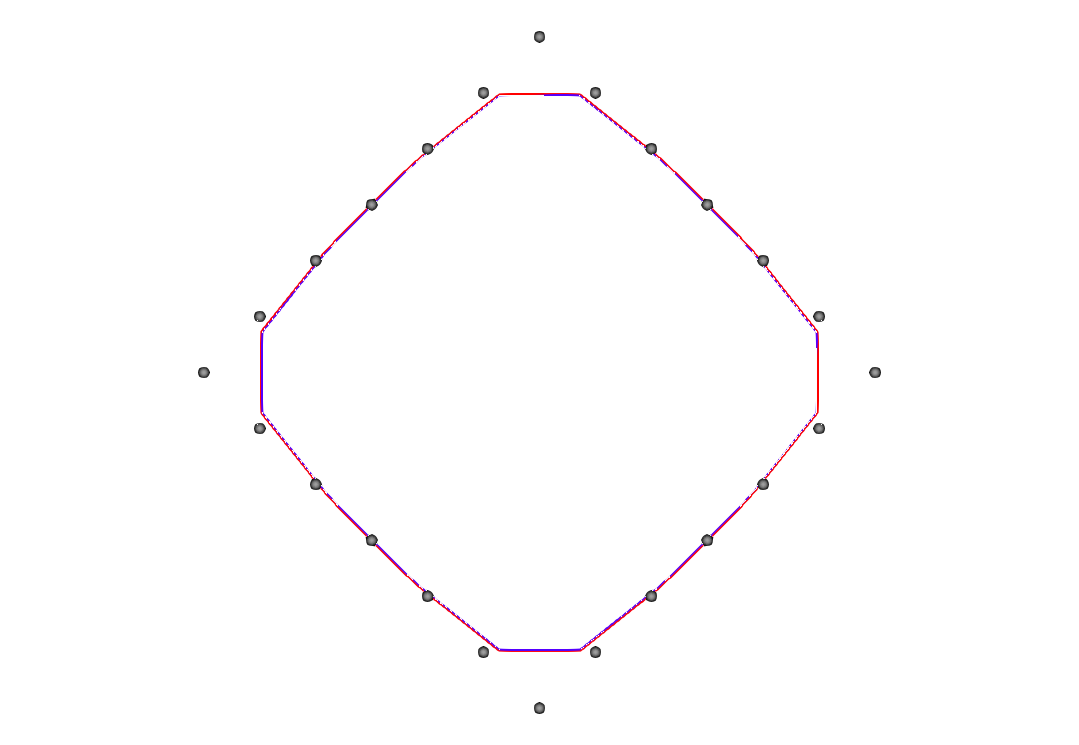}
\end{minipage}\hfill
\begin{minipage}{0.33\linewidth}
    \centering
    \includegraphics[width=1.0\linewidth]{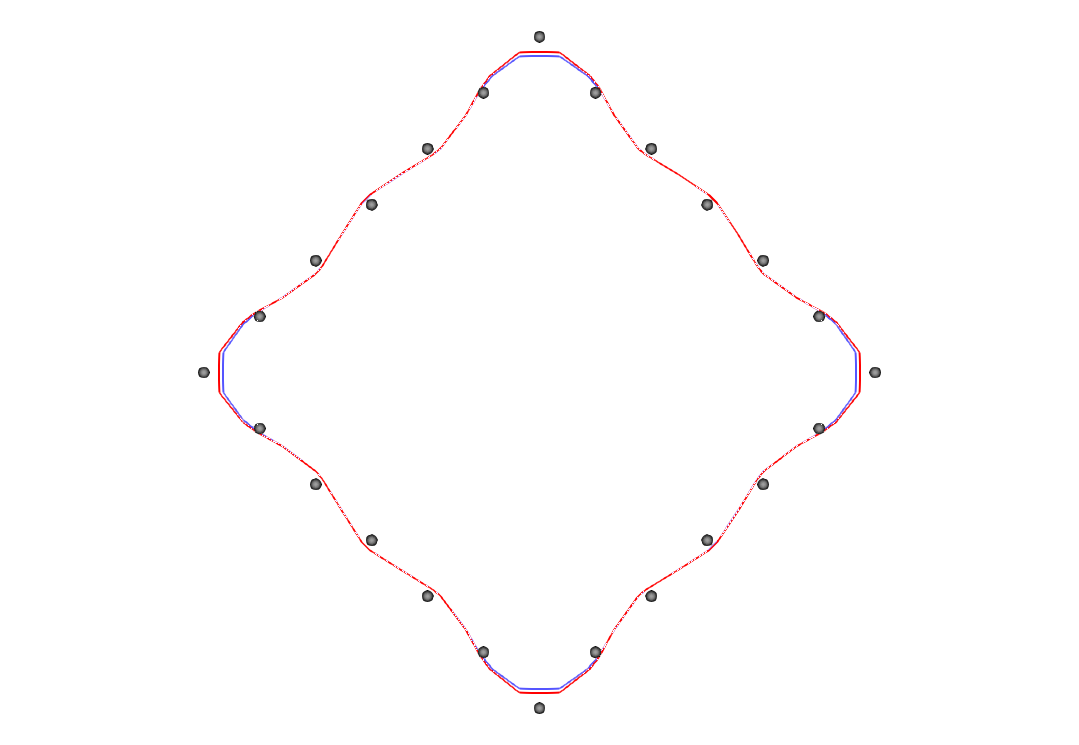}
\end{minipage}\hfill
\begin{minipage}{0.33\linewidth}
    \centering
    \includegraphics[width=1.0\linewidth]{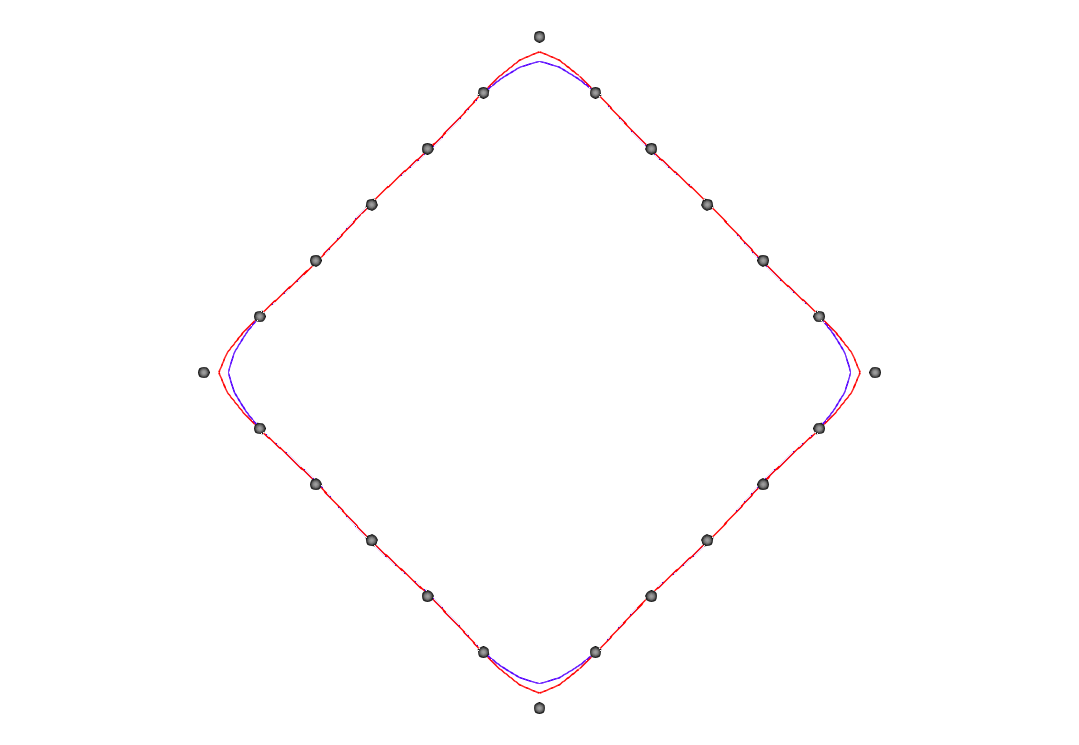}
\end{minipage}\\
\caption{From left to right: a comparison between the curves reconstructed with $\Qone$ (blue line) and $\WENO$ (red line) interpolators at each run.}
\label{fig:rhombusReconstruction}
\end{figure}

\begin{figure}
\begin{minipage}{0.45\linewidth}
    \begin{tikzpicture}

\begin{semilogyaxis}[%
  width=4.5cm,height=4cm,scale only axis,
  xtick={0, 43, 117, 143},
  xticklabels={$0$, $43$, $117$, $143$},
  major tick length = 0,
  yticklabel style={/pgf/number format/sci},
  ymin=1e-1,xmin=0,
  legend columns=-1, 
  legend style={
  	legend pos=north east,
  	draw=none,
  },
  title={$E_2(\phi^n)$}
]
\addplot[color=blue] table[col sep=space]{grafici/squareEn1.txt};
\addplot[color=red] table[col sep=space]{grafici/squareEn2.txt};
\addplot[color=green] table[col sep=space]{grafici/squareEn3.txt};

\draw[dashed,help lines] (axis cs:43,1e-20) -- (axis cs:43,100) ;
\draw[dashed,help lines] (axis cs:117,1e-20) -- (axis cs:117,100) ;
\draw[dashed,help lines] (axis cs:143,1e-20) -- (axis cs:143,100) ;
\draw[stealth-stealth,help lines] (axis cs:1,2e-1) -- node[pos=0.5,below,black]{\tiny $r=1$} (axis cs:42,2e-1) ;
\draw[stealth-stealth,help lines] (axis cs:44,2e-1) -- node[pos=0.5,below,black]{\tiny $r=2$} (axis cs:116,2e-1) ;
\draw[stealth-stealth,help lines] (axis cs:118,2e-1) -- node[pos=0.5,below,black]{\tiny $r=3$} (axis cs:142,2e-1) ;

\end{semilogyaxis}
\end{tikzpicture}
\end{minipage}\hspace{0.8cm}
\begin{minipage}{0.45\linewidth}
    \begin{tikzpicture}

\begin{semilogyaxis}[%
  width=4.5cm,height=4cm,scale only axis,
  xtick={0, 43, 117, 143},
  xticklabels={$0$, $43$, $117$, $143$},
  major tick length = 0,
  yticklabel style={/pgf/number format/sci},
  ymin=5e-3,
  xmin=0,
  legend columns=-1, 
  legend style={
  	legend pos=north east,
  	draw=none,
  },
  title={$Err_{\pcloud}^n$}
]
\addplot[color=blue] table[col sep=space]{grafici/squareErrCloud1.txt};
\addplot[color=red] table[col sep=space]{grafici/squareErrCloud2.txt};
\addplot[color=green] table[col sep=space]{grafici/squareErrCloud3.txt};

\draw[dashed,help lines] (axis cs:43,1e-20) -- (axis cs:43,100) ;
\draw[dashed,help lines] (axis cs:117,1e-20) -- (axis cs:117,100) ;
\draw[dashed,help lines] (axis cs:143,1e-20) -- (axis cs:143,100) ;
\draw[stealth-stealth,help lines] (axis cs:1,1e-2) -- node[pos=0.5,below,black]{\tiny $r=1$} (axis cs:42,1e-2) ;
\draw[stealth-stealth,help lines] (axis cs:44,1e-2) -- node[pos=0.5,below,black]{\tiny $r=2$} (axis cs:116,1e-2) ;
\draw[stealth-stealth,help lines] (axis cs:118,1e-2) -- node[pos=0.5,below,black]{\tiny $r=3$} (axis cs:142,1e-2) ;

\end{semilogyaxis}
\end{tikzpicture}
\end{minipage}\\

    \caption{Energy and error on cloud computed for the square test using $\WENO$ interpolator.}
    \label{fig:rhombusGraphs}
\end{figure}
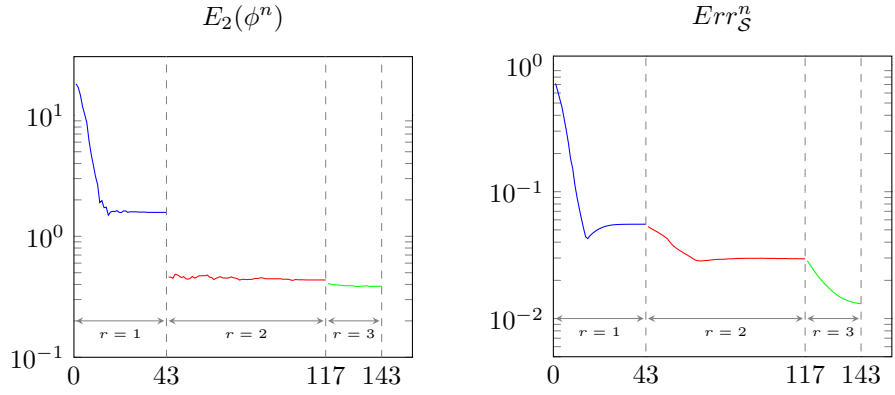

\begin{table}
\begin{center}
\tiny
\begin{tabular}{|c|c|c|c|c|c|c|c|}
\hline & & \multicolumn{3}{|c|}{$\Qone$} &\multicolumn{3}{|c|}{$\WENO$} \\
\hline $r$ & Grid & $L^{1}$-err & Energy & Error on $\pcloud$ & $L^{1}$-err & Energy & Error on $\pcloud$ \\ \hline \hline
$1$ & $30\times30$ & $1.58e+00$ & $5.09e-01$ & $6.65e-02$ & $1.51e+00$ & $4.57e-01$ & $5.54e-02$\\ \hline 
$2$ & $34\times34$ & $1.66e-01$ & $4.28e-01$ & $3.14e-02$ & $1.55e-01$ & $4.36e-01$ & $2.96e-02$\\ \hline 
$3$ & $51\times51$ & $6.52e-02$ & $3.97e-01$ & $2.13e-02$ & $4.00e-02$ & $3.88e-01$ & $1.31e-02$\\ \hline 
$4$ & $85\times85$ & $2.76e-02$ & $4.06e-01$ & $2.08e-02$ & $1.35e-02$ & $4.00e-01$ & $1.32e-02$\\ \hline 
$5$ & $153\times153$ & $1.21e-02$ & $4.12e-01$ & $1.65e-02$ & $5.53e-03$ & $4.07e-01$ & $1.23e-02$\\ \hline 
\end{tabular}
\end{center}
\caption{Errors and energy computed for the square test at the end of each run. To better appreciate the results, Algorithm~\ref{algo} has been performed five times.}
\label{tab:errGlobSquare}
\end{table}

In the second test we consider a cloud made up of $24$ points representing a square rotated by 45 degrees with respect to the Cartesian axis. The reconstruction takes respectively $43$, $74$, $26$ iterations and $2.78e-01$ seconds in the $\Qone$ case, and $39$, $97$, $18$ iterations and $2.71e-01$ seconds in the $\WENO$ case. The final curves of each run are depicted in Fig.~\ref{fig:rhombusReconstruction}. Graphs of the energy functional and of the error on the cloud are reported in Fig.~\ref{fig:rhombusGraphs}, exclusively for the $\WENO$ case. It is worth highlighting the role of $\delta$: the curvature regularization in the last run penalizes the sharp corners, thus providing a final contour with a controlled maximum curvature, which, although the higher resolution of the grid, is not passing exactly through the corresponding points of the cloud.
Table~\ref{tab:errGlobSquare} shows the errors computed with respect to the exact data. We point out that in this test the $\WENO$ based algorithm has the same computational time of the $\Qone$ one but leads to lower errors. The main differences are near the corners, as it can be seen from the right panel of Fig.~\ref{fig:rhombusReconstruction}.

\subsection{Synthetic 3d data sets}
In this subsection we present numerical tests for the reconstruction of 3d shapes. Here the data sets are synthetic and are made up by sampling points on simple geometrical forms for which the exact signed distance function is known.

\subsubsection{Sphere}
\begin{figure}
\begin{minipage}{0.33\linewidth}
    \centering
    \includegraphics[width=1.0\linewidth]{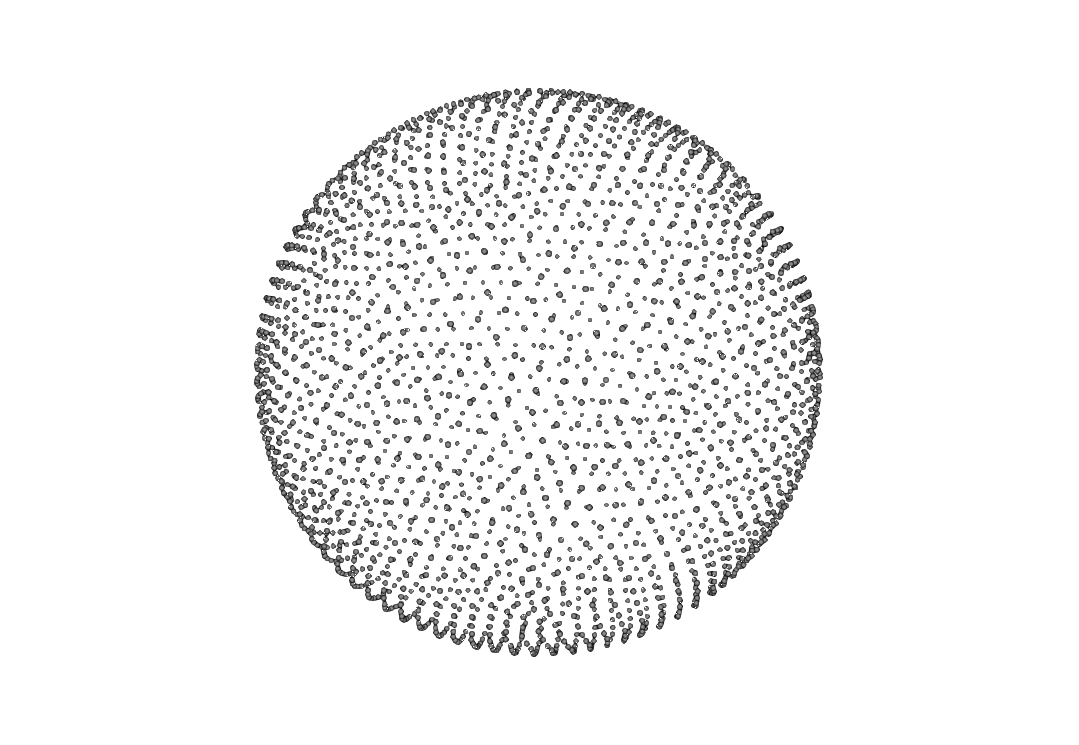}
\end{minipage}\hfill
\begin{minipage}{0.33\linewidth}
    \centering
    \includegraphics[width=1.0\linewidth]{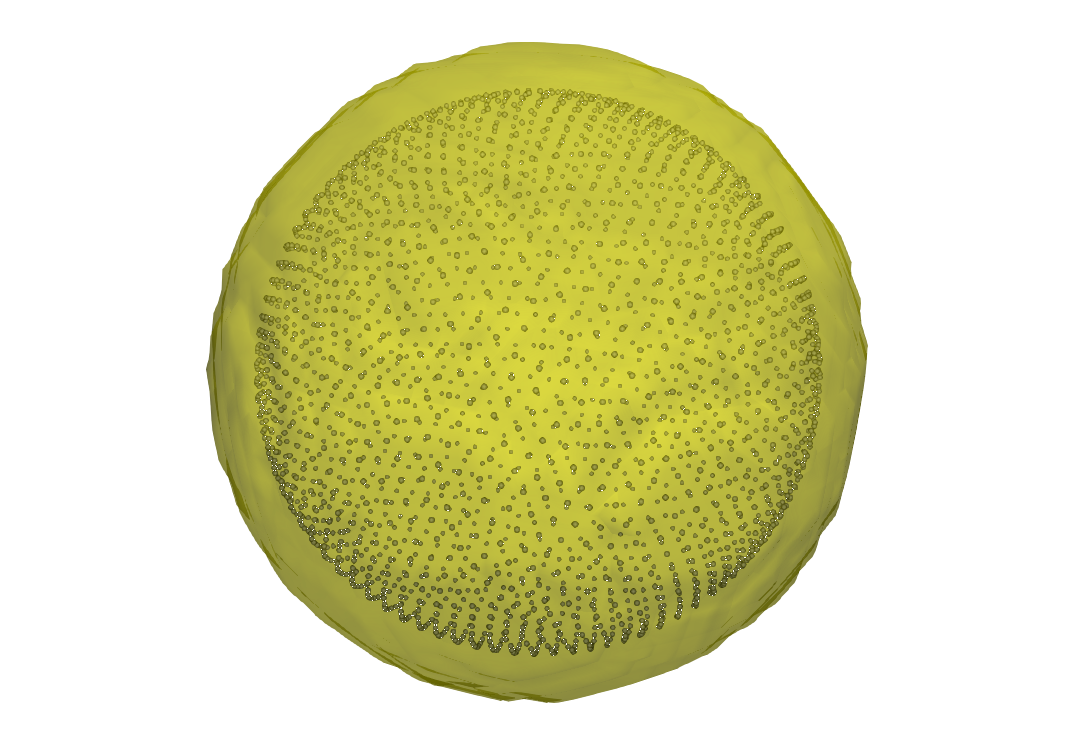}
\end{minipage}\hfill
\begin{minipage}{0.33\linewidth}
    \centering
    \includegraphics[width=1.0\linewidth]{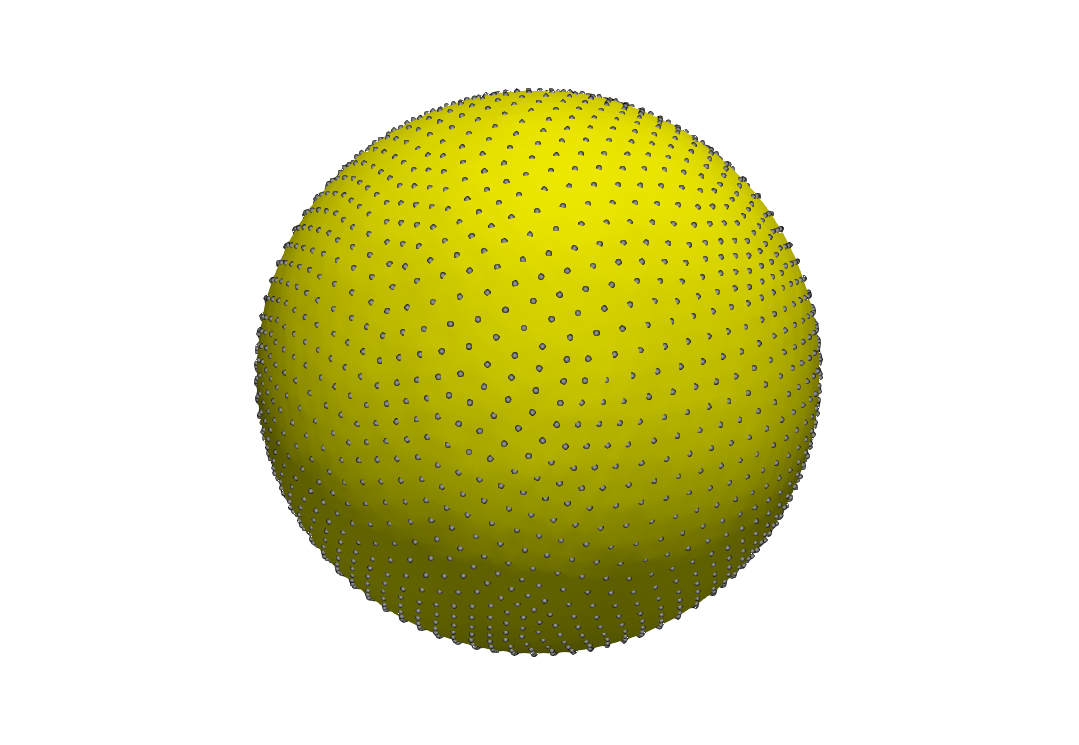}
\end{minipage}\\
\begin{minipage}{0.33\linewidth}
    \centering
    \includegraphics[width=1.0\linewidth]{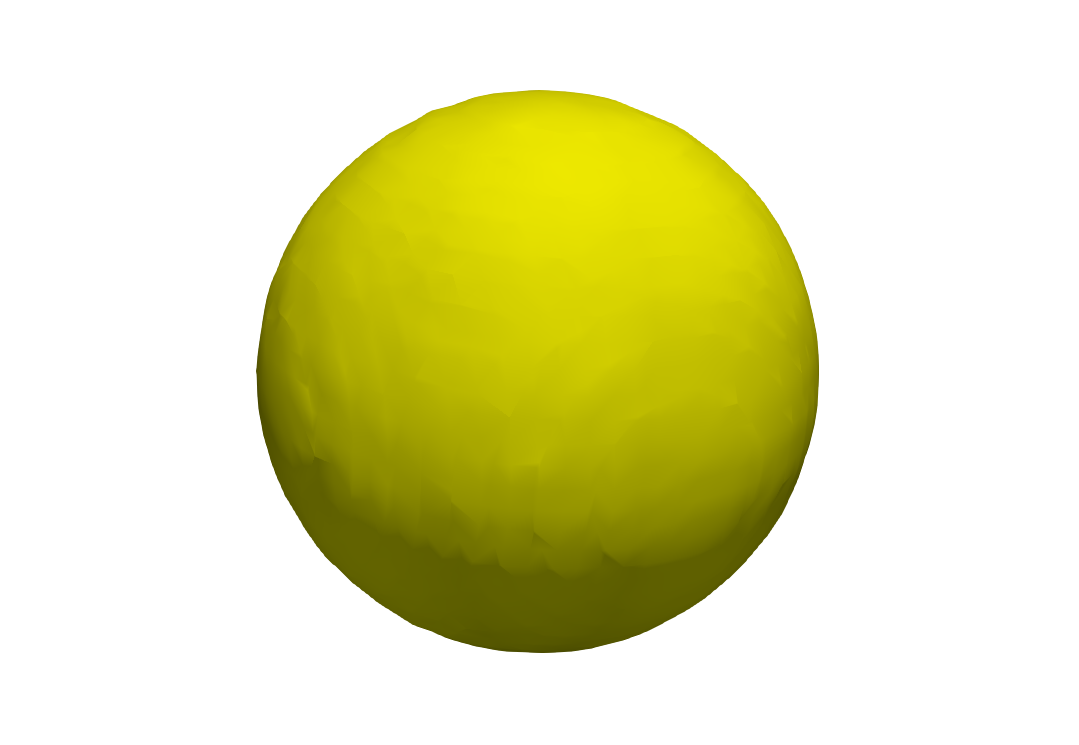}
\end{minipage}\hfill
\begin{minipage}{0.33\linewidth}
    \centering
    \includegraphics[width=1.0\linewidth]{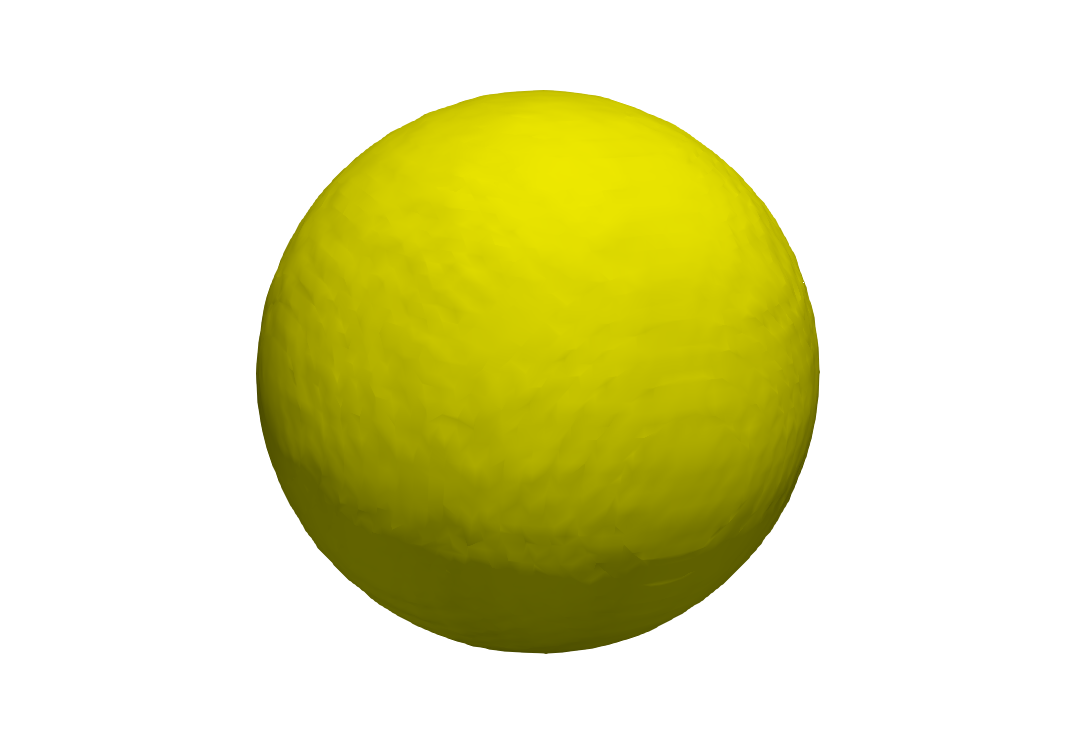}
\end{minipage}\hfill
\begin{minipage}{0.33\linewidth}
    \centering
    \includegraphics[width=1.0\linewidth]{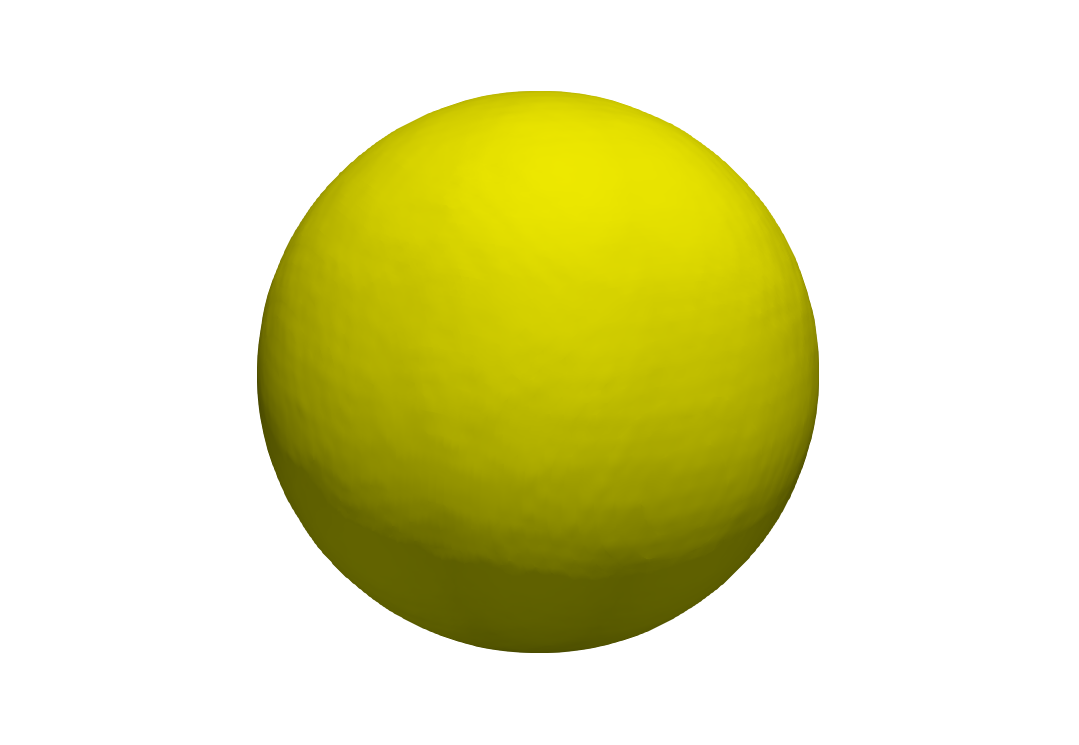}
\end{minipage}
\caption{Reconstruction steps of a sphere for the $\WENO$ case. In the first row the point cloud, the initial surface and the final one are depicted. In the second row the final surfaces obtained in each progressive run are depicted without the data set.}
\label{fig:sfera}
\end{figure}

The first 3d test has been performed on a point cloud made up of $2562$ points chosen on a sphere of radius $1$. The reconstruction procedure is illustrated in Fig.~\ref{fig:sfera} for the $\WENO$ case. In this simple case it is worth emphasizing that the role of the curvature regularization is to provide a smoother final surface than the ones provided by the previous two runs. Looking a Table~\ref{tab:errGlobSphere} one can compare the results obtained with different interpolators and observe the better performance of the $\WENO$ interpolator. 

\begin{table}
\begin{center}
\tiny
\begin{tabular}{|c|c|c|c|c|c|c|}
\hline & \multicolumn{3}{|c|}{$\Qone$} &\multicolumn{3}{|c|}{$\WENO$} \\
\hline $r$ & $L^{1}$-err & Energy & Error on $\pcloud$ & $L^{1}$-err & Energy & Error on $\pcloud$ \\ \hline \hline
$1$ & $9.14e-02$ & $2.17e-01$ & $3.66e-03$ & $8.76e-02$ & $1.94e-01$ & $3.09e-03$\\ \hline 
$2$ & $1.93e-02$ & $1.72e-01$ & $2.64e-03$ & $1.72e-02$ & $1.69e-01$ & $2.33e-03$\\ \hline 
$3$ & $8.35e-03$ & $1.61e-01$ & $2.67e-03$ & $6.74e-03$ & $1.56e-01$ & $2.05e-03$\\ \hline 
\end{tabular}
\end{center}
\caption{Errors and energy computed for the sphere at the end of each run.}
\label{tab:errGlobSphere}
\end{table}

\subsubsection{Cube\&Spheres}
The second test concerns a three-dimensional domain composed by a cube joined with three spheres. The cube edge length was $0.8$, the first sphere has radius $0.25$ and centre at the middle of an edge of the cube, while the other two had radius $0.15$ and were centred onto the two vertices of the opposite edge of the cube. The geometrical object was rotated in such a way that no face nor edge were aligned with the background Cartesian grid.

\begin{figure}
\centering
\begin{minipage}{0.45\linewidth}
    \centering
    \includegraphics[width=1.0\linewidth]{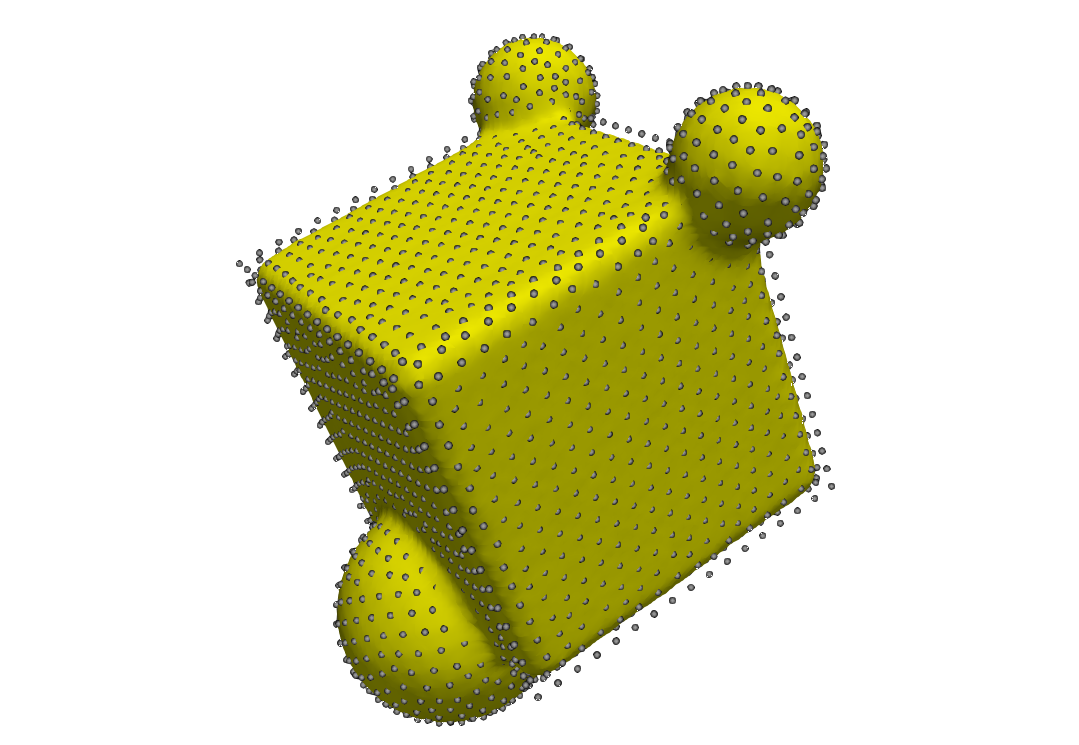}
\end{minipage}
\begin{minipage}{0.45\linewidth}
    \centering
    \includegraphics[width=1.0\linewidth]{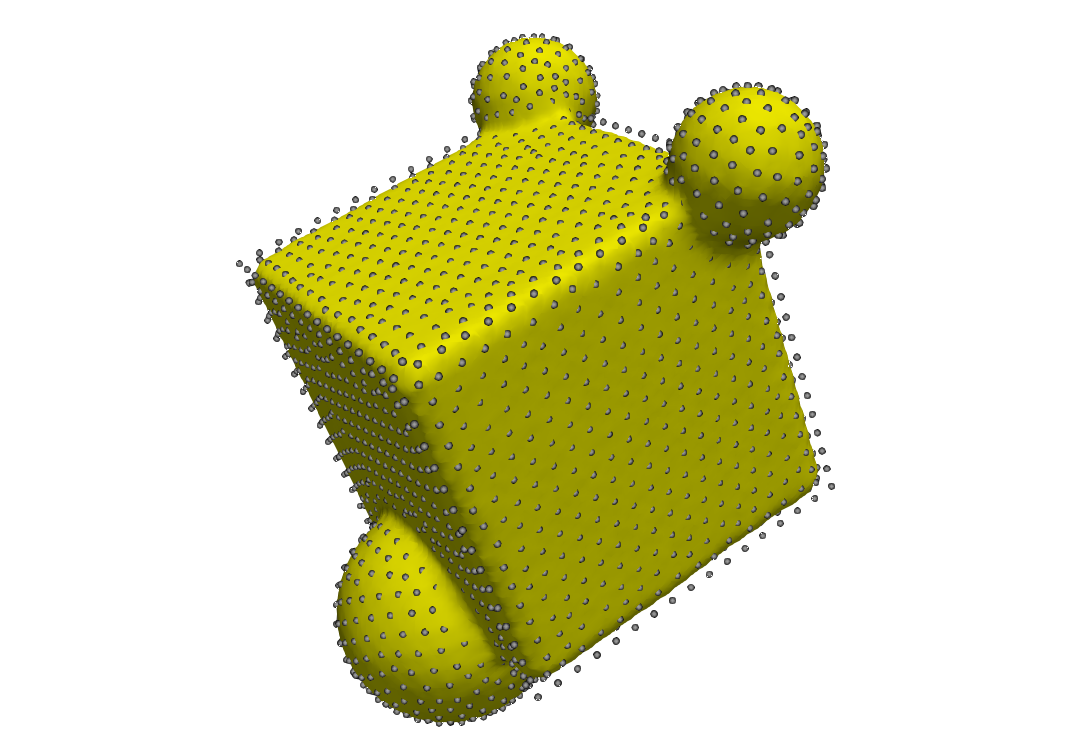}
\end{minipage}
\caption{The reconstructed surfaces of the ``Cube\&Spheres''. On the left, the result obtained with the Q1 interpolant, while on the right, the result obtained using the $\WENO$ interpolant.}
\label{fig:recCubosfere}
\end{figure}

\begin{table}
\begin{center}
\tiny
\begin{tabular}{|c|c|c|c|c|c|c|}
\hline & \multicolumn{3}{|c|}{$\Qone$} &\multicolumn{3}{|c|}{$\WENO$} \\
\hline $r$ & $L^{1}$-err & Energy & Error on $\pcloud$ & $L^{1}$-err & Energy & Error on $\pcloud$ \\ \hline \hline
$1$ & $2.78e-02$ & $7.79e-02$ & $5.60e-03$ & $2.68e-02$ & $7.14e-02$ & $4.74e-03$\\ \hline 
$2$ & $4.61e-03$ & $6.46e-02$ & $2.79e-03$ & $4.52e-03$ & $6.40e-02$ & $2.65e-03$\\ \hline 
$3$ & $1.90e-03$ & $6.11e-02$ & $3.50e-03$ & $1.66e-03$ & $5.98e-02$ & $2.74e-03$\\ \hline  
\end{tabular}
\end{center}
\caption{Errors and energy computed for the Cube\&Spheres at the end of each run.}
\label{tab:errGlobCubosfere}
\end{table}

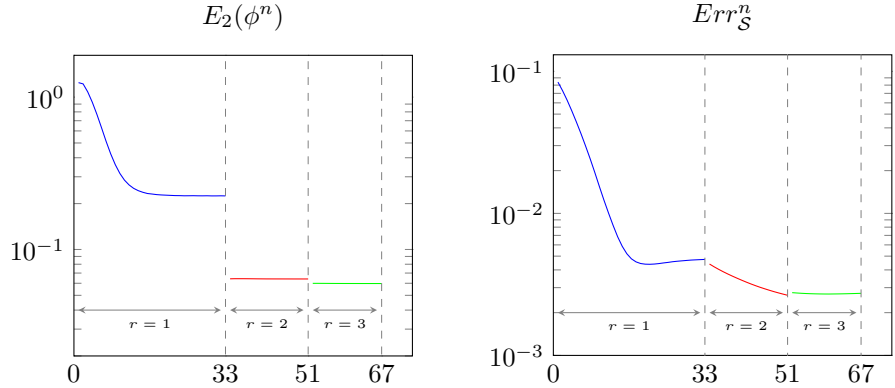
\begin{figure}
\begin{minipage}{0.45\linewidth}
    \begin{tikzpicture}

\begin{semilogyaxis}[%
  width=4.5cm,height=4cm,scale only axis,
  xtick={0, 33, 51, 67},
  xticklabels={$0$, $33$, $51$, $67$},
  major tick length = 0,
  yticklabel style={/pgf/number format/sci},
  ymin=2e-2,
  xmin=0,
  legend columns=-1, 
  legend style={
  	legend pos=north east,
  	draw=none,
  },
  title={$E_2(\phi^n)$}
]
\addplot[color=blue] table[col sep=space]{grafici/cubosfereEn1.txt};
\addplot[color=red] table[col sep=space]{grafici/cubosfereEn2.txt};
\addplot[color=green] table[col sep=space]{grafici/cubosfereEn3.txt};

\draw[dashed,help lines] (axis cs:33,1e-20) -- (axis cs:33,100) ;
\draw[dashed,help lines] (axis cs:51,1e-20) -- (axis cs:51,100) ;
\draw[dashed,help lines] (axis cs:67,1e-20) -- (axis cs:67,100) ;
\draw[stealth-stealth,help lines] (axis cs:1,4e-2) -- node[pos=0.5,below,black]{\tiny $r=1$} (axis cs:32,4e-2) ;
\draw[stealth-stealth,help lines] (axis cs:34,4e-2) -- node[pos=0.5,below,black]{\tiny $r=2$} (axis cs:50,4e-2) ;
\draw[stealth-stealth,help lines] (axis cs:52,4e-2) -- node[pos=0.5,below,black]{\tiny $r=3$} (axis cs:66,4e-2) ;

\end{semilogyaxis}
\end{tikzpicture}
\end{minipage}\hspace{0.8cm}
\begin{minipage}{0.45\linewidth}
    \begin{tikzpicture}

\begin{semilogyaxis}[%
  width=4.5cm,height=4cm,scale only axis,
  xtick={0, 33, 51, 67},
  xticklabels={$0$, $33$, $51$, $67$},
  major tick length = 0,
  yticklabel style={/pgf/number format/sci},
  ymin=1e-3,
  xmin=0,
  legend columns=-1, 
  legend style={
  	legend pos=north east,
  	draw=none,
  },
  title={$Err_{\pcloud}^n$}
]
\addplot[color=blue] table[col sep=space]{grafici/cubosfereErrCloud1.txt};
\addplot[color=red] table[col sep=space]{grafici/cubosfereErrCloud2.txt};
\addplot[color=green] table[col sep=space]{grafici/cubosfereErrCloud3.txt};

\draw[dashed,help lines] (axis cs:33,1e-20) -- (axis cs:33,100) ;
\draw[dashed,help lines] (axis cs:51,1e-20) -- (axis cs:51,100) ;
\draw[dashed,help lines] (axis cs:67,1e-20) -- (axis cs:67,100) ;
\draw[stealth-stealth,help lines] (axis cs:1,2e-3) -- node[pos=0.5,below,black]{\tiny $r=1$} (axis cs:32,2e-3) ;
\draw[stealth-stealth,help lines] (axis cs:34,2e-3) -- node[pos=0.5,below,black]{\tiny $r=2$} (axis cs:50,2e-3) ;
\draw[stealth-stealth,help lines] (axis cs:52,2e-3) -- node[pos=0.5,below,black]{\tiny $r=3$} (axis cs:66,2e-3) ;

\end{semilogyaxis}
\end{tikzpicture}
\end{minipage}\\
\caption{Energy and error on cloud computed for the ``Cube\&Spheres'' using $\WENO$ interpolator.}
\label{fig:err:cubosfere}
\end{figure}

Results are shown in Fig.~\ref{fig:recCubosfere} and in Table~\ref{tab:errGlobCubosfere}. As for the 2d case of the square, the curvature regularization provides rounded edges and corners in the final reconstruction, especially in the $\Qone$ case. Here too the $\WENO$ reconstruction gives lower errors according to all three measures and especially on the error on the point cloud $Err_{\pcloud}$. Plot of the energy functional and of the error on the cloud are shown in Fig.~\ref{fig:err:cubosfere} for the $\WENO$ case. Note in particular that the minimum of the error on $\pcloud$ is reached at the end of the second run and is about $2.65e-03$, while at the end of the third run the error is about $2.74e-3$, which is an effect of the curvature regularization.

\subsection{Data sets from laser scans}

Here we test our method using data sets coming from laser-scanning of real objects. These point clouds are made available in the Digital Shape WorkBench of the AIM@SHAPE and VISIONAIR  projects \cite{frog} or in the 3D Scanning Repository of the Stanford University \cite{standford}.

\subsubsection{The ``Frog'' and the ``Bunny''}

\begin{figure}
\centering
\begin{minipage}{0.45\linewidth}
    \centering
    \includegraphics[width=1.0\linewidth]{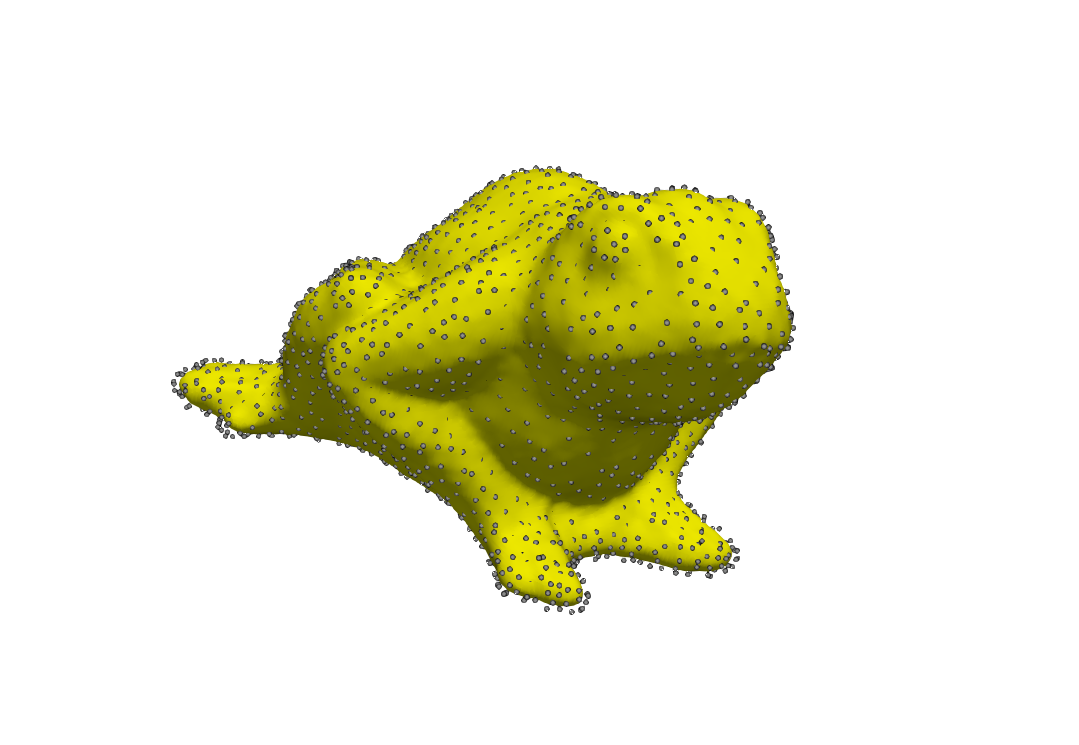}
\end{minipage}
\begin{minipage}{0.45\linewidth}
    \centering
    \includegraphics[width=1.0\linewidth]{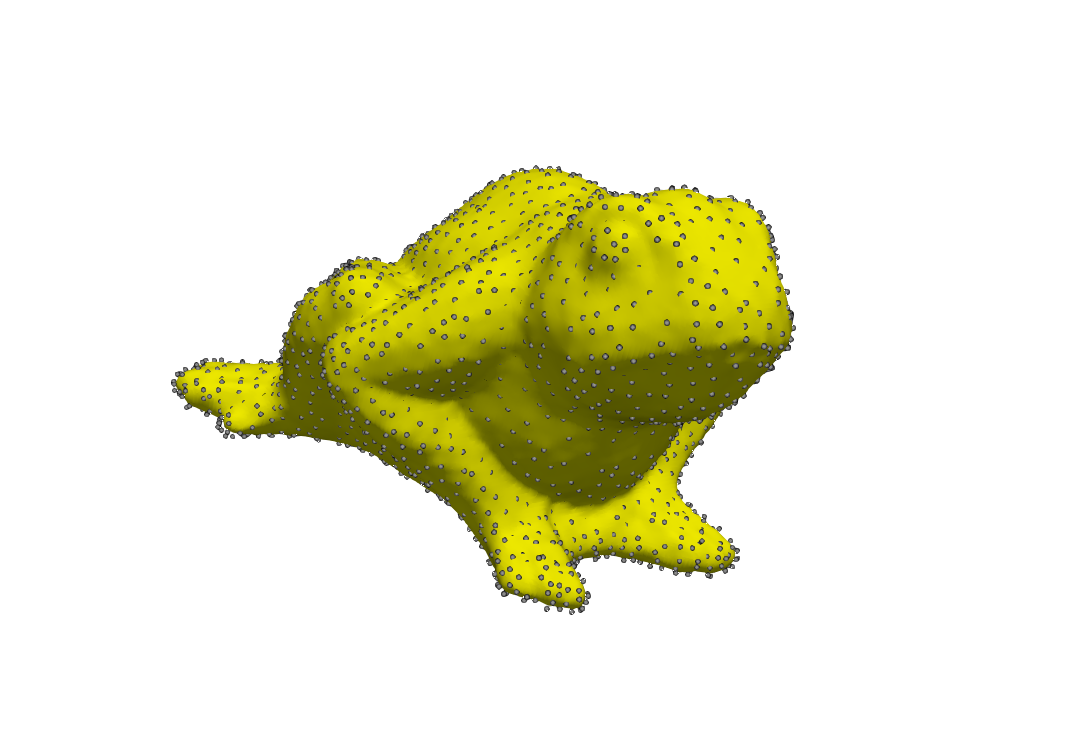}
\end{minipage}\\
\centering
\begin{minipage}{0.45\linewidth}
    \centering
    \includegraphics[width=1.0\linewidth]{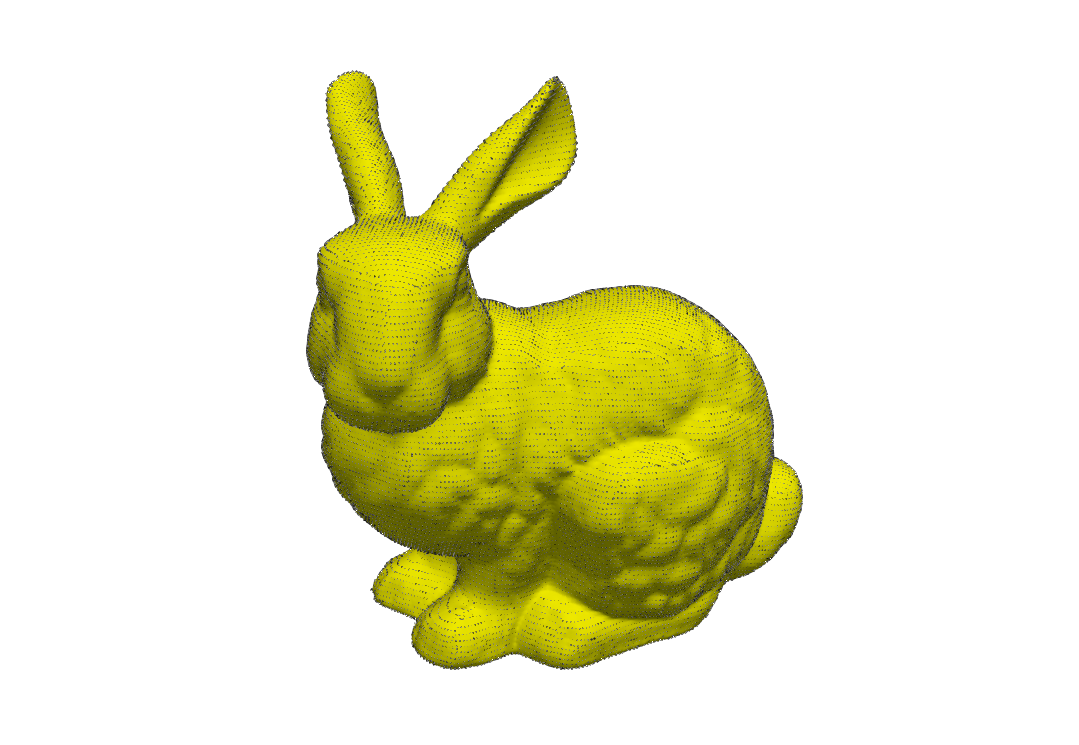}
\end{minipage}
\begin{minipage}{0.45\linewidth}
    \centering
    \includegraphics[width=1.0\linewidth]{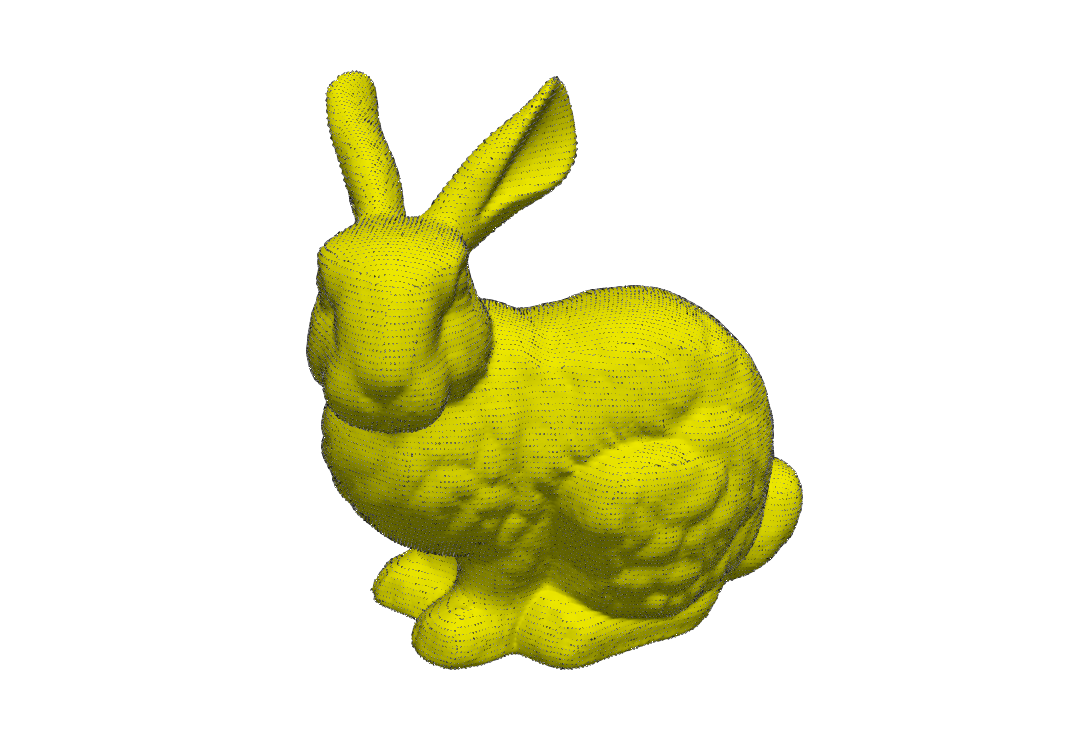}
\end{minipage}
\caption{The reconstructed surfaces of the ``Frog" and the ``Bunny", respectively on the first and on the second row. On the left, the result obtained with the $\Qone$ interpolant, while on the right, the result obtained using the $\WENO$ interpolant. Especially for the ``Frog'', due to the lower resolution of the point cloud, one can appreciate how the WENO recostruction recovers the details of the paws.}
\label{fig:recFrogBunny}
\end{figure}

\begin{table}
\begin{center}
\tiny
\begin{tabular}{|c|c|c|c|c|}
\hline & \multicolumn{2}{|c|}{Q1} &\multicolumn{2}{|c|}{WENO} \\
\hline $r$ & Energy & Error on $\pcloud$ & Energy & Error on $\pcloud$ \\ \hline \hline
$1$ & $1.41e-01$ & $6.39e-03$ & $1.35e-01$ & $5.58e-03$\\ \hline 
$2$ & $1.31e-01$ & $3.76e-03$ & $1.28e-01$ & $3.02e-03$\\ \hline 
$3$ & $1.32e-01$ & $4.31e-03$ & $1.31e-01$ & $3.08e-03$\\ \hline 
\end{tabular}
\end{center}
\caption{Errors and energy computed for the ``Frog'' at the end of each run.}
\label{tab:errGlobFrog}
\end{table}

\begin{table}
\begin{center}
\tiny
\begin{tabular}{|c|c|c|c|c|}
\hline & \multicolumn{2}{|c|}{Q1} &\multicolumn{2}{|c|}{WENO} \\
\hline $r$ & Energy & Error on $\pcloud$ & Energy & Error on $\pcloud$ \\ \hline \hline
$1$ & $8.68e-03$ & $3.07e-04$ & $8.36e-03$ & $2.95e-04$\\ \hline 
$2$ & $7.98e-03$ & $1.35e-04$ & $7.86e-03$ & $1.41e-04$\\ \hline 
$3$ & $7.56e-03$ & $1.47e-04$ & $7.45e-03$ & $1.22e-04$\\ \hline 
\end{tabular}
\end{center}
\caption{Errors and energy computed for the ``Bunny'' at the end of each run.}
\label{tab:errGlobBunny}
\end{table}

We first consider a point cloud named ``Frog'' \cite{frog} and one named ``Bunny'' \cite{standford}. The results are shown in Fig.~\ref{fig:recFrogBunny} and in Tables~\ref{tab:errGlobFrog},~\ref{tab:errGlobBunny}. In both these two cases the data sets have holes in the bottom and that's why we had to increase the threshold $\gamma_{\pcloud}$ for the initial data computation, thus we set $K_{\pcloud} = 5.0$ for the ``Frog'', and $K_{\pcloud} = 10.0$ for the ``Bunny''. This obviously produced an enlarged initial contour and therefore more iterations were needed during the first run, but this simple strategy secured us from getting distorted results deriving from a bad initial data. Of course, alternatively, one could also fill these holes in an ad-hoc preprocessing stage by adding some points to the cloud.

We also stress as the curvature term in the last run has produced a final level set reconstruction that is everywhere smooth, but still retains many details of the shape: see for example the legs and feet of both shapes and the ears of the bunny.

\subsubsection{A teapot with tiny details}

\begin{figure}
\centering
\begin{minipage}{0.45\linewidth}
    \centering
    \includegraphics[width=1.0\linewidth]{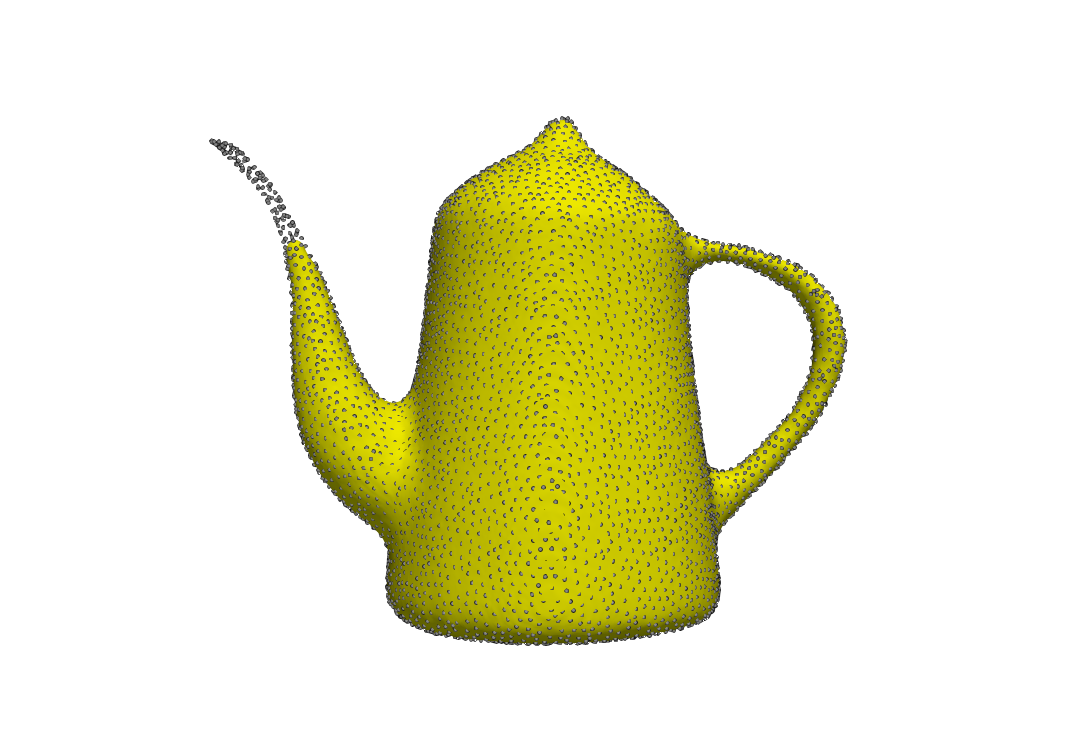}
\end{minipage}
\begin{minipage}{0.45\linewidth}
    \centering
    \includegraphics[width=1.0\linewidth]{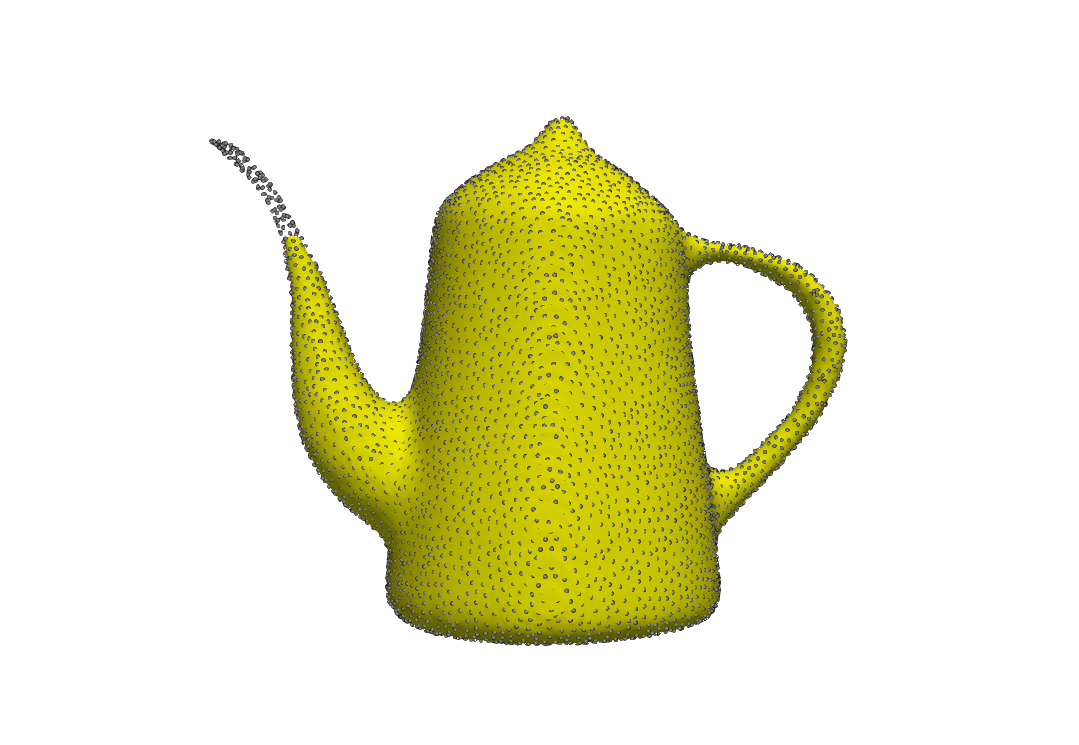}
\end{minipage}\\
\begin{minipage}{0.45\linewidth}
    \centering
    \includegraphics[width=1.0\linewidth]{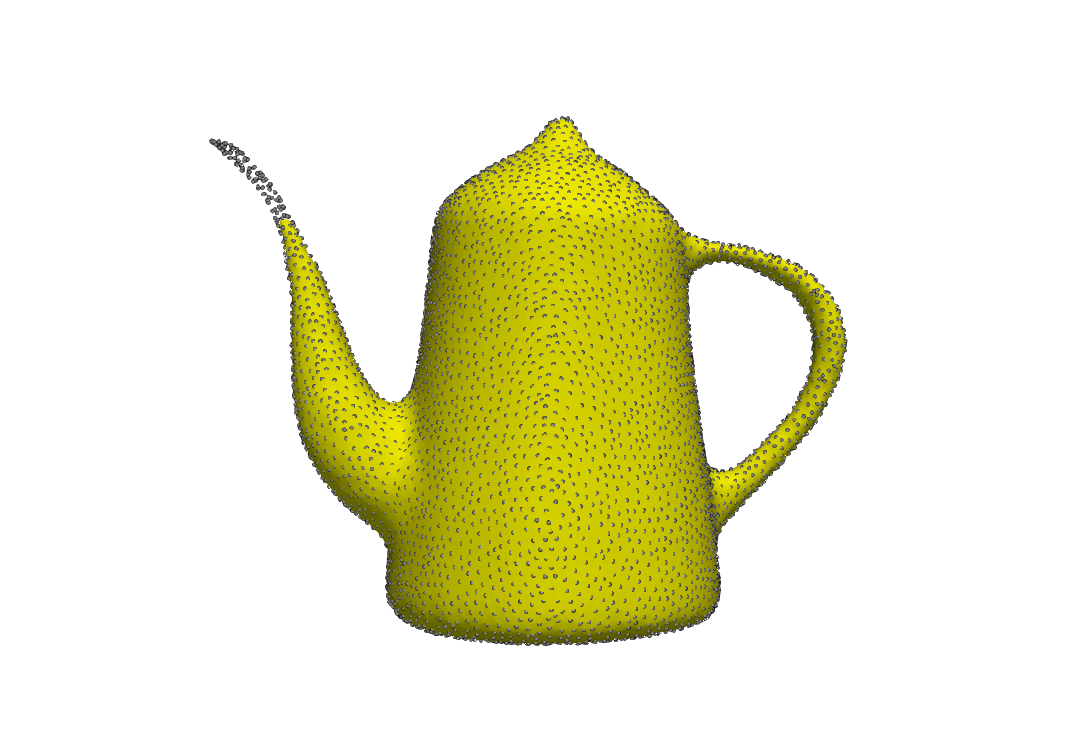}
\end{minipage}
\begin{minipage}{0.45\linewidth}
    \centering
    \includegraphics[width=1.0\linewidth]{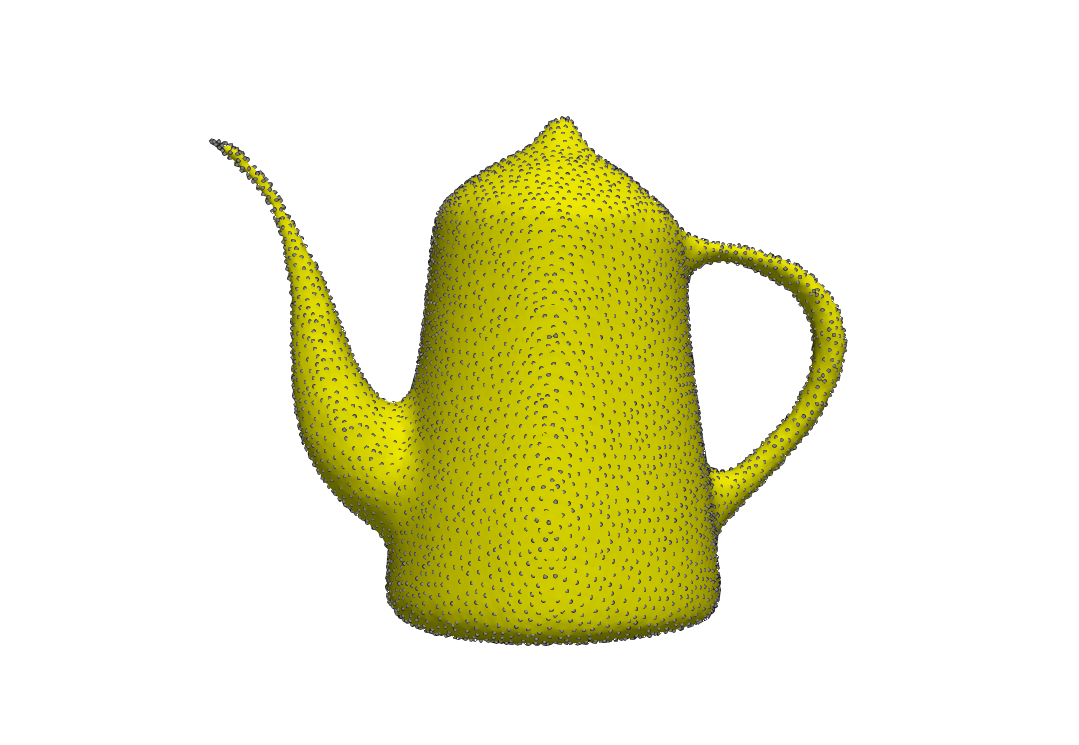}
\end{minipage}
\caption{Final reconstructions of the ``Teapot''. Left panels collect results obtained with $\Qone$ interpolator, while on the right $\WENO$ interpolator is used. In the first row the usual setting of the parameters \eqref{eq:dxdt} fails in recovering the shape of the spout, while in the second row one can appreciate that, starting with a smaller grid, and then increasing the resolution as usual, reconstructions are better, especially the $\WENO$ one that is able to recover the detail of the spout.}
\label{fig:teapot}
\end{figure}

In the next test we focus our attention on the ability of our algorithm to capture tiny details of an object. A point cloud named ``Teapot'' \cite{frog} has been considered for this aim since it presents a more complex topology and parts whose cross section is more or less comparable with the size of the cloud. In the first row of Fig.~\ref{fig:teapot}, we can see that with the usual setting of the parameters \eqref{eq:dxdt} we miss the spout of the teapot. This can be at least partially fixed by starting the algorithm with a smaller grid size in the first run, setting $\dx^{(r)}=\frac{1}{2^{r}}h_{\pcloud}$. The second row of Fig.~\ref{fig:teapot} reports the results obtained with this last setting. Note that, however, it is necessary to involve the $\WENO$ reconstruction to have a full recovery of the details.
This better result of course requires a greater computational effort: the bottom row reconstruction takes about four times as long as the one above (see Table~\ref{tab:3dTimes}).

\subsubsection{Complex shapes and topology: the ``Happy Buddha'' and the ``Dragon''}

\begin{figure}
\centering
\begin{minipage}{0.3\linewidth}
\hspace*{-1.5cm}
    \includegraphics[width=2.2\linewidth]{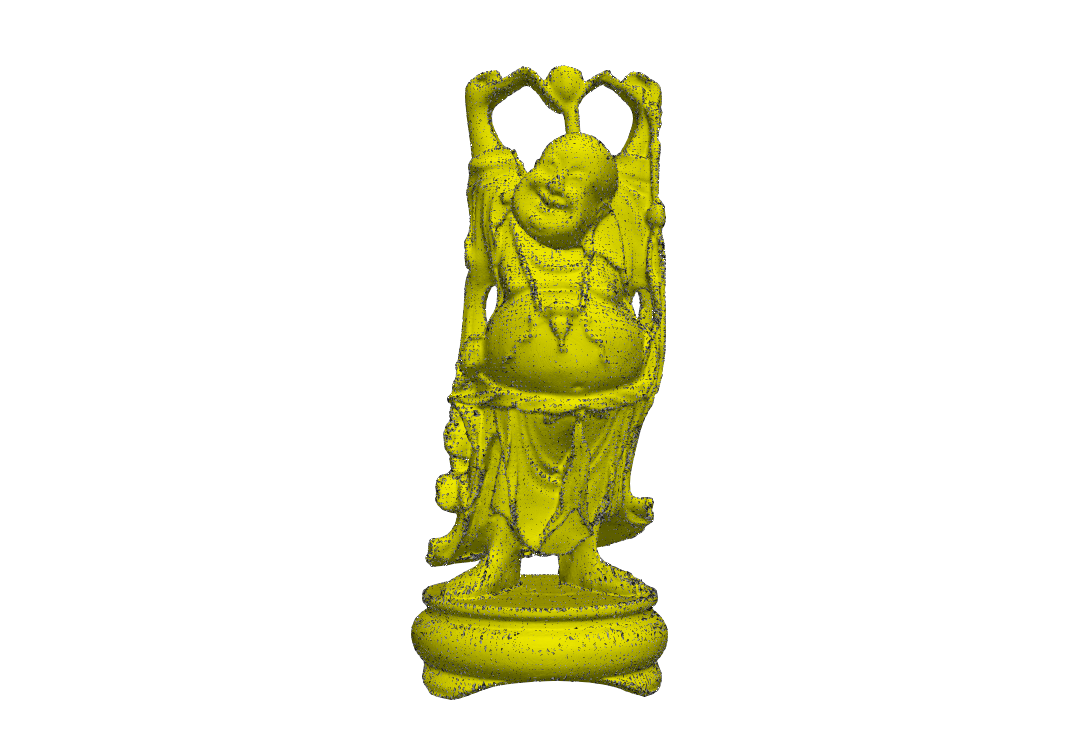}
\end{minipage}
\begin{minipage}{0.65\linewidth}
    \centering
    \includegraphics[width=0.9\linewidth]{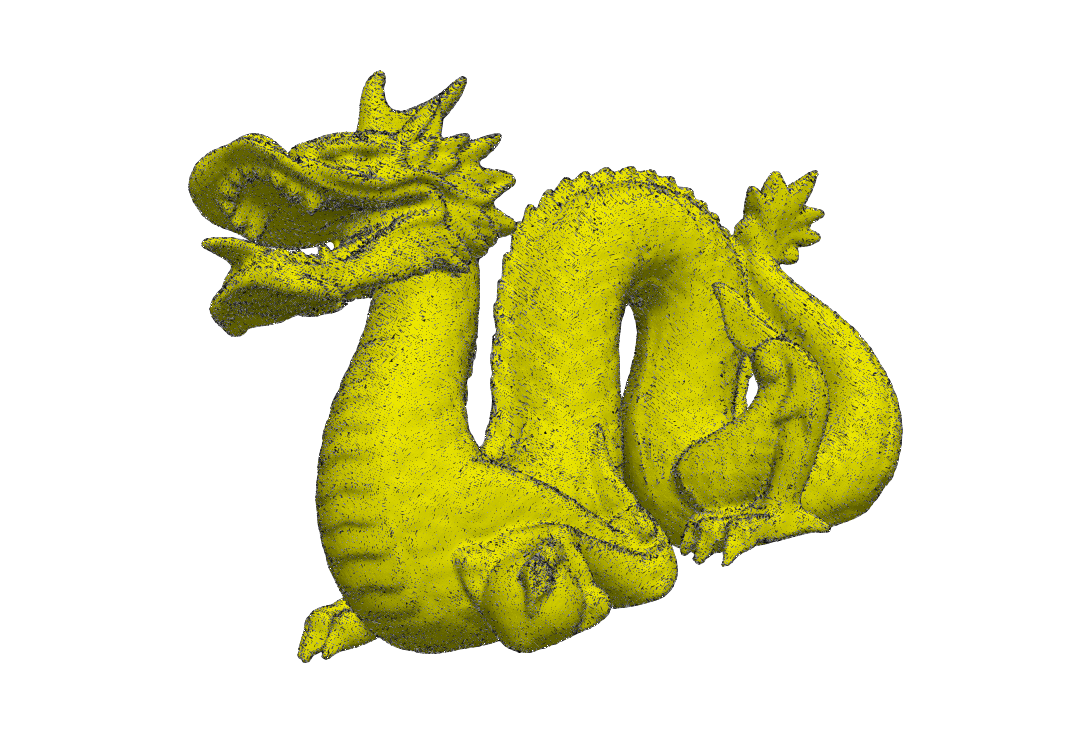}
\end{minipage}\\
    \caption{Final reconstructions of the ``Happy Buddha'' and the ``Dragon'' using $\WENO$ interpolator.}
    \label{fig:happy&dragon}
\end{figure}

Finally we have done some tests on very complex point clouds present in the Stanford 3D Scanning Repository \cite{standford}. We considered two data sets, named ``Happy Buddha'' and ``Dragon'' respectively, because they present some nice features like free holes, small bridges due to the carving and many details to be recovered. The results of our reconstructions are depicted in Fig.~\ref{fig:happy&dragon} and has been obtained setting $K_{\pcloud} = 10.0$ in the computation of the initial data, and $\dx^{(r)}=\frac{1}{2^{r-2}}h_{\pcloud}$ in order to limit the huge amount of grid points.

On the ``Happy Buddha'' test, we point out that the reconstructed surface recovered equally well both the very flat surfaces of the base and the tiny details on the sides, on the belly and on the drapery. Also all the holes of the highly nontrivial topology were correctly recovered, including the small ones on the sides: we point out that the initial data had detected the two big holes at the top but missed the two tiny ones at the sides which were recovered during the surface evolution, confirming the ability of the method to deal with topological changes of the surface.

On the ``Dragon'' test, we point out that the scales on the skin are still well approximated despite the curvature regularization. Also, we stress how the level set has been pushed inside the mouth during its evolution and how the sharp teeth shapes have been well approximated.

\subsection{Scalability of the algorithm}

Our algorithm, thanks to the localization of the computational effort on a band around the evolving zero level set, which is a codimension $1$ variety of $\R^n$, has a cost that scales as $\Ogrande(1/\dx^{n-1})$ under grid refinement. This is confirmed by comparing the two ``Teapot'' experiments (Table~\ref{tab:3dGrid} and~\ref{tab:3dTimes}). The memory footprint of our implementation scales instead as $\Ogrande(1/\dx^n)$, due to the data attached to the full grid $\grid$. This has suggested to consider a distributed memory parallel implementation based on the MPI paradigm.

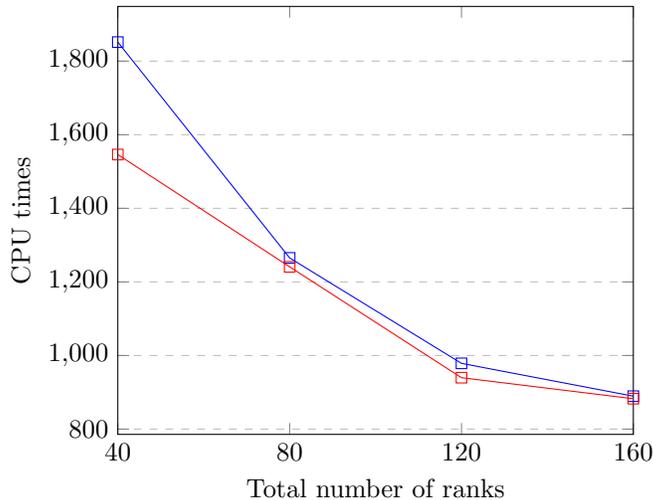
\begin{figure}
\centering
    \begin{tikzpicture}
\begin{axis}[
    xlabel={Total number of ranks},
    ylabel={CPU times},
    xmin=40, xmax=160,
    xtick={0,40,80,120,160},
    legend pos=north west,
    ymajorgrids=true,
    grid style=dashed,
]

\addplot[
    color=blue,
    mark=square,
    ]
    coordinates {
    (40,1.8517e+03)(80,1.2652e+03)(120,9.7864e+02)(160,8.8949e+02)
    };

\addplot[
    color=red,
    mark=square,
    ]
    coordinates {
    (40,1.5461e+03)(80,1.2410e+03)(120,9.3970e+02)(160,8.8258e+02)
    };

\end{axis}
\end{tikzpicture}
    \caption{Scalability results for the ``Dragon'' reconstruction. The blue marks report the CPU times in seconds required using a fixed number of ranks per node ($40$) and gradually increasing the number of nodes, from $1$ to $4$. The red marks refer to the opposite test: the number of nodes is fixed at $4$ and the number of ranks per node goes from $10$ to $40$.}
\label{fig:scalability}
\end{figure}

In this last subsection we have collected some scalability results to evaluate the efficiency of the algorithm. The test depicted in Fig.~\ref{fig:scalability} has been performed using the ``Dragon'' point cloud progressively increasing the number of nodes, while keeping constant the number of processors per node, and the other way round fixing the number of nodes and increasing the number of processors per node. 

Computing the parallel efficiency of a run with $M$ cores as $\epsilon = \frac{40 \cdot T_{40}}{M \cdot T_{M}}$, we get results between $60-75\%$ for the $80$ cores run and around $50\%$ for the ones involving more tasks. This can be ascribed to the interplay between the Cartesian grid partitioning for $\grid$ and the localization techniques of \S\ref{ssec:mask}. In fact, when the number of cores is increased, also the chance of having cores with no (or small) intersection with the active computational band is also increased, with a negative impact on load balancing.

\begin{table}
\begin{center}
\tiny
\begin{tabular}{|c|c|c|c|c|c|c|c|c|}
\hline
\multirow{2}{*}{Point cloud} & \multicolumn{3}{c|}{Grid} \\
\cline{2-4}
& $r=1$ & $r=2$ & $r=3$ \\
\hline
Sphere & $54\times54\times54$ & $78\times78\times78$ & $135\times135\times135$\\\hline
Cube\&Spheres & $55\times62\times55$ & $81\times94\times81$ & $141\times166\times141$\\\hline
Frog & $79\times65\times77$ & $116\times89\times113$ & $211\times157\times204$\\\hline
Bunny & $197\times195\times162$ & $332\times329\times262$ & $642\times637\times503$\\\hline
Teapot & $103\times65\times86$ & $176\times101\times143$ & $330\times181\times264$\\\hline
Teapot (finer) & $184\times109\times151$ & $330\times181\times264$ & $638\times340\times506$\\\hline
Happy Buddha & $282\times641\times282$ & $522\times1241\times523$ & $1023\times2460\times1024$\\\hline
Dragon & $504\times364\times243$ & $966\times687\times444$ & $1910\times1353\times866$\\\hline
\end{tabular}
\end{center}
\caption{Dimensions of the Cartesian grid involved for each 3d test at each progressive run.}
\label{tab:3dGrid}
\end{table}

\begin{table}
\begin{center}
\begin{tabular}{|P{2cm}|P{1.5cm}|P{1.5cm}|}
\hline
\multirow{2}{*}{Point cloud} & \multicolumn{2}{c|}{Total CPU time} \\
\cline{2-3}
& $\Qone$ & $\WENO$ \\
\hline
Sphere & 0.09 & 0.26 \\\hline
Cube\&Spheres & 0.14 & 0.33 \\\hline
Frog & 0.31 & 0.88 \\\hline
Bunny & 1.78 & 3.42 \\\hline
Teapot & 0.24 & 0.50 \\\hline
Teapot (finer) & 0.95 & 1.78 \\\hline
Happy Buddha & 16.95 & 57.55 \\\hline
Dragon & 12.93 & 25.68 \\\hline
\end{tabular}
\end{center}
\caption{Computational times (min) of the algorithm.}
\label{tab:3dTimes}
\end{table}

\section{Conclusions}
\label{sec:conclusions}

We presented a complete workflow for the reconstruction of surfaces starting from unorganized point cloud data, without assuming any information about the connection between points. We used the level set method to detect and evolve these surfaces in an implicit way. In particular, we used the signed distance function as level set function to preserve numerical accuracy and to get some nice geometrical features for free. The numerical method is based on a semi-Lagrangian scheme with a local $\Qone$ interpolator or a $\WENO$ one and thus it works locally and on fixed rectangular grids. Many efforts have been made to save computational costs. The code runs in parallel and it turned out to be quite efficient despite the very fine resolution of the grids employed and the high amount of grid points involved in the computations. Numerical tests performed for the 2d cases and more complex 3d cases suggested that the workflow illustrated here is quite promising, especially using $\WENO$ techniques, which shows better results with respect to the simple $\Qone$ interpolator. Future research should be directed towards improving the efficiency of the algorithm involving quadree and octree adaptive meshes to solve the load balancing issues highlighted in \S$4.4$, focusing all the computational effort in the regions close to the data. Interesting directions of research would be also the replacement of the current semi-Lagrangian scheme with its monotone version (see \cite{CaFe:2008}) to achieve a monotone decrease of the error and of the energy functional.

\paragraph{Acknowledgements} The authors would like to thank the PETSc developer team for promptly including in the 3.19 release the relevant features of DMSWARM that made possible the simulations in this paper.

The authors acknowledge the CINECA award under the ISCRA initiative, for the availability of high-performance computing resources and support.

Both authors are members of the Gruppo Nazionale Calcolo Scientifico-Istituto Nazionale di Alta Matematica (GNCS-INdAM).

\paragraph{Data availability}
The datasets analysed during the current study are available in the repositories cited in the text.

\paragraph{Declarations}
The authors declare that they have no conflict of interest.

\bibliographystyle{spmpsci}
\bibliography{levelset.bib}

\begin{thebibliography}{10}
\providecommand{\url}[1]{{#1}}
\providecommand{\urlprefix}{URL }
\expandafter\ifx\csname urlstyle\endcsname\relax
  \providecommand{\doi}[1]{DOI~\discretionary{}{}{}#1}\else
  \providecommand{\doi}{DOI~\discretionary{}{}{}\begingroup
  \urlstyle{rm}\Url}\fi

\bibitem{frog}
AIM@SHAPE, VISIONAIR: Visionair shape repository.
\newblock \url{http://visionair.ge.imati.cnr.it/ontologies/shapes/} (2022).
\newblock
  \urlprefix\url{http://visionair.ge.imati.cnr.it/ontologies/shapes/view.jsp?id=268-frog_-_merged}.
\newblock Accessed March 2022

\bibitem{Am:delaunay2004}
Amenta, N., Bern, M., Kellis, M.: A new voronoi-based surface reconstruction
  algorithm.
\newblock Proc. Siggraph '98 \textbf{98} (2004).
\newblock \doi{10.1145/280814.280947}

\bibitem{petsc-user-ref}
Balay, S., Abhyankar, S., Adams, M.F., Brown, J., Brune, P., Buschelman, K.,
  Dalcin, L., Eijkhout, V., Gropp, W.D., Karpeyev, D., Kaushik, D., Knepley,
  M.G., May, D.A., McInnes, L.C., Mills, R.T., Munson, T., Rupp, K., Sanan, P.,
  Smith, B.F., Zampini, S., Zhang, H., Zhang, H.: {PETSc/TAO} users manual.
\newblock Tech. Rep. ANL-21/39 - Revision 3.19, Argonne National Laboratory
  (2023)

\bibitem{petsc-efficient}
Balay, S., Gropp, W.D., McInnes, L.C., Smith, B.F.: Efficient management of
  parallelism in object oriented numerical software libraries.
\newblock In: E.~Arge, A.M. Bruaset, H.P. Langtangen (eds.) Modern Software
  Tools in Scientific Computing, pp. 163--202. Birkh{\"{a}}user Press (1997)

\bibitem{BeTaSe:survey2016}
Berger, M., Tagliasacchi, A., Seversky, L., Alliez, P., Guennebaud, G., Levine,
  J., Sharf, A., Silva, C.: A survey of surface reconstruction from point
  clouds.
\newblock Computer Graphics Forum \textbf{36}, n/a--n/a (2016).
\newblock \doi{10.1111/cgf.12802}

\bibitem{BeTaSe:survey2014}
Berger, M., Tagliasacchi, A., Seversky, L., Alliez, P., Levine, J., Sharf, A.,
  Silva, C.: State of the art in surface reconstruction from point clouds.
\newblock Eurographics 2014-State of the Art Reports  (2014)

\bibitem{BCCF:21}
Bonaventura, L., Calzola, E., Carlini, E., Ferretti, R.: Second order fully
  semi-{L}agrangian discretizations of advection-diffusion-reaction systems.
\newblock J. Sci. Comput. \textbf{88}(23) (2021).
\newblock \doi{10.1007/s10915-021-01518-8}

\bibitem{BF:14:SLdiffusion}
Bonaventura, L., Ferretti, R.: Semi-{L}agrangian methods for parabolic problems
  in divergence form.
\newblock {SIAM} J. Sci. Comput. \textbf{36}(5), A2458--A2477 (2014).
\newblock \doi{10.1137/140969713}

\bibitem{CaFaFe:2010}
Carlini, E., Falcone, M., Ferretti, R.: Convergence of a large time-step scheme
  for mean curvature motion.
\newblock Interfaces and Free Boundaries \textbf{12}(4), 409--411 (2010).
\newblock \doi{10.4171/IFB/240}

\bibitem{CaFe:2008}
Carlini, E., Ferretti, R.: A semi-lagrangian approximation of min–max type
  for the stationary mean curvature equation.
\newblock In: Numerical Mathematics and Advanced Applications, pp. 679--686
  (2008).
\newblock \doi{10.1007/978-3-540-69777-0_81}

\bibitem{CaFe:2017}
Carlini, E., Ferretti, R.: A semi-lagrangian scheme with radial basis
  approximation for surface reconstruction.
\newblock Comput. Vis. Sci. \textbf{18}(2-3), 103--112 (2017).
\newblock \doi{10.1007/s00791-016-0274-2}

\bibitem{CFR05}
Carlini, E., Ferretti, R., Russo, G.: A weighted essentially nonoscillatory,
  large time-step scheme for {H}amilton-{J}acobi equations.
\newblock {SIAM} J. Sci. Comput. \textbf{27}(3), 1071 -- 1091 (2006).
\newblock \doi{10.1137/040608787}

\bibitem{RBF1}
Carr, J., Beatson, R., Cherrie, J., Mitchell, T., Fright, W., Mccallum, B.,
  Evans, T.: Reconstruction and representation of 3d objects with radial basis
  functions.
\newblock ACM SIGGRAPH  (2001).
\newblock \doi{10.1145/383259.383266}

\bibitem{RBF2}
Carr, J., Beatson, R., McCallum, B., Fright, W., McLennan, T., Mitchell, T.:
  Smooth surface reconstruction from noisy range data.
\newblock In: Proceedings of the 1st International Conference on Computer
  Graphics and Interactive Techniques in Australasia and South East Asia,
  GRAPHITE '03, pp. 119–126+297 (2003).
\newblock \doi{10.1145/604471.604495}

\bibitem{CaFr:surveyDelaunay}
Cazals, F., Giesen, J.: Delaunay Triangulation Based Surface Reconstruction,
  pp. 231--276 (2006).
\newblock \doi{10.1007/978-3-540-33259-6_6}

\bibitem{leastSquare2}
Cheng, Z.Q., Wang, Y.Z., Li, B., Xu, K., Dang, G., Jin, S.Y.: A survey of
  methods for moving least squares surfaces.
\newblock pp. 9--23 (2008).
\newblock \doi{10.2312/VG/VG-PBG08/009-023}

\bibitem{cdss:mach19}
Coco, A., Preda, S., Semplice, M.: From point clouds to 3d simulations of
  marble sulfation.
\newblock p. 153 – 174 (2023).
\newblock \doi{10.1007/978-981-99-3679-3_10}

\bibitem{CR2018}
Coco, A., Russo, G.: Second order finite-difference ghost-point multigrid
  methods for elliptic problems with discontinuous coefficients on an arbitrary
  interface.
\newblock J. Comput. Phys. \textbf{361}, 299--330 (2018).
\newblock \doi{https://doi.org/10.1016/j.jcp.2018.01.016}

\bibitem{CSS:monum}
Coco, A., Semplice, M., Serra~Capizzano, S.: A level-set multigrid technique
  for nonlinear diffusion in the numerical simulation of marble degradation
  under chemical pollutants.
\newblock Appl. Math. \& Comput. \textbf{386}, 125503 (2020).
\newblock \doi{10.1016/j.amc.2020.125503}

\bibitem{DaMi:surfMesh2015}
Daniel, P., Medl’a, M., Mikula, K., Remesikova, M.: Reconstruction of
  surfaces from point clouds using a {L}agrangian surface evolution model.
\newblock pp. 589--600 (2015).
\newblock \doi{10.1007/978-3-319-18461-6_47}

\bibitem{Ed:delaunay1998}
Edelsbrunner, H.: Shape reconstruction with {D}elaunay complex.
\newblock pp. 119--132 (1998).
\newblock \doi{10.1007/BFb0054315}

\bibitem{FaFe:2003}
Falcone, M., Ferretti, R.: Consistency of a large time-step scheme for mean
  curvature motion.
\newblock In: F.~Brezzi, A.~Buffa, S.~Corsaro, A.~Murli (eds.) Numerical
  Mathematics and Advanced Applications, pp. 495--502. Springer Milan (2003).
\newblock \doi{10.1007/978-88-470-2089-4_46}

\bibitem{GF:2002}
Gibou, F., Fedkiw, R.P., Cheng, L.T., Kang, M.: A second-order-accurate
  symmetric discretization of the {P}oisson equation on irregular domains.
\newblock J. Comput. Phys. \textbf{176}(1), 205--227 (2002).
\newblock \doi{https://doi.org/10.1006/jcph.2001.6977}

\bibitem{NURBS2}
Gregorski, B., Hamann, B., Joy, K.: Reconstruction of b-spline surfaces from
  scattered data points.
\newblock In: Proceedings Computer Graphics International 2000, pp. 163--170
  (2000).
\newblock \doi{10.1109/CGI.2000.852331}

\bibitem{HaMeSc:2008}
Hartmann, D., Meinke, M., Schröder, W.: Differential equation based
  constrained reinitialization for level set methods.
\newblock Journal of Computational Physics \textbf{227}(14), 6821--6845 (2008).
\newblock \doi{https://doi.org/10.1016/j.jcp.2008.03.040}.
\newblock
  \urlprefix\url{https://www.sciencedirect.com/science/article/pii/S0021999108001964}

\bibitem{He:2019}
He, Y., Huska, M., Kang, S.H., Liu, H.: Fast algorithms for surface
  reconstruction from point cloud.
\newblock In: Springer Proceedings in Mathematics and Statistics, vol. 360, pp.
  61--80 (2021).
\newblock \doi{10.1007/978-981-16-2701-9_4}

\bibitem{hoppe1992}
Hoppe, H., Derose, T., Duchamp, T., Mcdonald, J., Stuet-zle, W.: Surface
  reconstruction from unorganized point clouds.
\newblock Computer Graphics pp. 71--78 (1992).
\newblock \doi{https://doi.org/10.1145/142920.134011}

\bibitem{JiangShu:96}
Jiang, G.S., Shu, C.W.: Efficient implementation of weighted {ENO} schemes.
\newblock J. Comput. Phys. \textbf{126}, 202--228 (1996).
\newblock \doi{10.1006/jcph.1996.0130}

\bibitem{Kosa2017}
Kósa, B., Haličková-Brehovská, J., Mikula, K.: New efficient numerical
  method for 3d point cloud surface reconstruction by using level set methods.
\newblock Proceedings of Equadiff 2017 Conference pp. 387--396 (2017).
\newblock
  \urlprefix\url{http://www.iam.fmph.uniba.sk/amuc/ojs/index.php/equadiff/article/view/798}

\bibitem{standford}
Laboratory, S.U.C.G.: The stanford 3d scanning repository.
\newblock \urlprefix\url{http://graphics.stanford.edu/data/3Dscanrep/}.
\newblock Accessed March 2022

\bibitem{DL:2023apr}
Liu, X.: Research on 3{D} object reconstruction method based on deep learning.
\newblock Highlights in Science, Engineering and Technology \textbf{39},
  1221--1227 (2023).
\newblock \doi{10.54097/hset.v39i.6732}

\bibitem{LSApp:2004}
Marcon, M., Piccarreta, L., Sarti, A., Tubaro, S.: Fast point-cloud wrapping
  through level-set evolution.
\newblock In: 1st European Conference on Visual Media Production (CVMP) 2004,
  pp. 119--125 (2004)

\bibitem{OshStFe:2004}
Osher, S., Fedkiw, R.: The Level Set Methods and Dynamic Implicit Surfaces,
  vol.~57, pp. xiv+273 (2004).
\newblock \doi{10.1115/1.1760520}

\bibitem{OshSet:1988}
Osher, S., Sethian, J.A.: Fronts propagating with curvature-dependent speed:
  Algorithms based on hamilton-jacobi formulations.
\newblock J. Computat. Phys. \textbf{79}(1), 12--49 (1988).
\newblock \doi{10.1016/0021-9991(88)90002-2}

\bibitem{NURBS1}
Park, I.K., Yun, I.D., Lee, S.U.: Constructing nurbs surface model from
  scattered and unorganized range data.
\newblock In: Second International Conference on 3-D Digital Imaging and
  Modeling (Cat. No.PR00062), pp. 312--320 (1999).
\newblock \doi{10.1109/IM.1999.805361}

\bibitem{PENG1999}
Peng, D., Merriman, B., Osher, S., Zhao, H., Kang, M.: A pde-based fast local
  level set method.
\newblock J. Computat. Phys. \textbf{155}(2), 410--438 (1999).
\newblock \doi{https://doi.org/10.1006/jcph.1999.6345}

\bibitem{Re:2011}
Remondino, F.: Heritage recording and 3d modeling with photogrammetry and 3d
  scanning.
\newblock Remote Sensing \textbf{3}, 1104--1138 (2011).
\newblock \doi{10.3390/rs3061104}

\bibitem{ReEl:2006}
Remondino, F., El-Hakim, S.: Image‐based 3d modelling: A review.
\newblock The Photogrammetric Record \textbf{21}, 269 -- 291 (2006).
\newblock \doi{10.1111/j.1477-9730.2006.00383.x}

\bibitem{Se:1999}
Sethian, J.: Level set methods and fast marching methods. evolving interfaces
  in computational geometry, fluid mechanics, computer vision, and materials
  science.
\newblock vol. 3. Cambridge university press, 1999. \textbf{3} (1999).
\newblock \doi{10.1090/S0025-5718-00-01345-4}

\bibitem{DL:2021}
Sharma, R., Schwandt, T., Kunert, C., Urban, S., Broll, W.: Point cloud
  upsampling and normal estimation using deep learning for robust surface
  reconstruction.
\newblock pp. 70--79 (2021).
\newblock \doi{10.5220/0010211600700079}

\bibitem{DL:2023gen}
Sulzer, R., Marlet, R., Vallet, B., Landrieu, L.: A survey and benchmark of
  automatic surface reconstruction from point clouds (2024).
\newblock \doi{https://doi.org/10.48550/arXiv.2301.13656}

\bibitem{SSO:1994}
Sussman, M., Smereka, P., Osher, S.: A level set approach for computing
  solutions to incompressible two-phase flow.
\newblock J. Computat. Phys. \textbf{114}(1), 146--159 (1994).
\newblock \doi{10.1006/jcph.1994.1155}

\bibitem{scan3d2020}
Wang, Q., Tan, Y., Mei, Z.: Computational methods of acquisition and processing
  of 3d point cloud data for construction applications.
\newblock Archives of Computational Methods in Engineering \textbf{27},
  479–499 (2020).
\newblock \doi{10.1007/s11831-019-09320-4}

\bibitem{RBF3}
Zeng, Y., Zhu, Y.: Implicit surface reconstruction based on a new interpolation
  / approximation radial basis function.
\newblock Computer Aided Geometric Design \textbf{92}, 102062 (2021).
\newblock \doi{10.1016/j.cagd.2021.102062}

\bibitem{Zhao2005AFS}
Zhao, H.: A fast sweeping method for eikonal equations.
\newblock Math. Comput. \textbf{74}, 603--627 (2005).
\newblock \doi{10.1090/S0025-5718-04-01678-3}

\bibitem{Zhao:2000}
Zhao, H.K., Osher, S., Merriman, B., Kang, M.: Implicit and nonparametric shape
  reconstruction from unorganized data using a variational level set method.
\newblock Computer Vision and Image Understanding \textbf{80}(3), 295--314
  (2000).
\newblock \doi{10.1006/cviu.2000.0875}

\end{thebibliography}
\end{document}